\documentclass[10pt]{article}
\usepackage{latexsym}
\usepackage{graphicx,subfigure,caption}
\usepackage{cite}
\usepackage{enumitem}
\usepackage{amssymb,mathrsfs,amsmath}
\usepackage{booktabs}
\usepackage[justification=centering]{caption}
\usepackage{subfigure}
\usepackage{float}
\usepackage{color}
\usepackage{cases}
\usepackage{subeqnarray}
\usepackage{cases}
\usepackage{endnotes}

\newtheorem{theorem}{Theorem}[part]
\newtheorem{definition}{Definition}[part]
\newtheorem{proposition}{Proposition}[part]

\newtheorem{remark}{Remark}[part]
\newtheorem{example}{Example}[part]

\def \ep{\hbox{ }\hfill$\Box$}

\begin{document}
\title{Block diagonalization of block circulant quaternion matrices and the fast calculation for T-product of quaternion tensors\thanks{This document is the results of the research project funded by National Natural Science Foundations of China (Grant No.11871472),  the Natural Science Foundation of Hunan Province, China (Grant No.2022JJ40543) and China Postdoctoral Science Foundation (Grant No.2021MD703978).} }

\author{Meng-Meng Zheng,\thanks{Department of Mathematics, National University of Defense Technology, Changsha, Hunan 410073, China. Email: mengmeng\_zheng@163.com.}
\ and\ Guyan Ni\thanks{Corresponding Author. Department of Mathematics, National University of Defense Technology, Changsha, Hunan 410073, China. Email: guyan-ni@163.com.}}
\date{}
\maketitle

\begin{abstract}
\noindent
With the great success of the T-product based real tensor methods in the color image and gray video processing, the establishment of T-product based quaternion tensor methods in the color video processing has encountered a challenge, which is the block diagonalization of block circulant quaternion matrices.
In this paper, we show that the discrete Fourier matrix $\mathbf{F_p}$ cannot diagonalize $p\times p$ circulant quaternion matrices, nor can the unitary quaternion matrices $\mathbf{F_p}\mathbf{j}$ and $\mathbf{F_p}(1+\mathbf{j})/\sqrt{2}$ with $\mathbf{j}$ being an imaginary unit of quaternion algebra, due to $1c = c1$ but $\mathbf{j}c = \bar{c}\mathbf{j}$ for any complex number $c$. Further, we establish sufficient and necessary conditions for a unitary quaternion matrix being a diagonalization matrix of circulant quaternion matrices, which shows that achieving the diagonalization of circulant quaternion matrices in the quaternion domain is too hard. Noting that  $\mathbf{l}q=\bar{q}\mathbf{l}$ holds for any quaternion $q$ with $\mathbf{l}$ being an imaginary unit of octonion algebra, we turn to the octonion domain for achieving the diagonalization of circulant quaternion matrices. We obtain several special multiplication cases between special octonions and quaternions. Based on that, the unitary octonion matrix $\mathbf{F_p}\mathbf{p}$ with $\mathbf{p}=\mathbf{l},\mathbf{il}$ or $(\mathbf{l}+\mathbf{il})/\sqrt{2}$ can diagonalize a circulant quaternion matrix of size $p\times p$, at the cost of $O(p\log p)$ via the fast Fourier transform (FFT); and unitary matrices $\mathbf{F_p}\mathbf{p}\otimes \mathbf{I_m}$ and $\mathbf{F_p}\mathbf{p}\otimes \mathbf{I_n}$ can block diagonalize a block circulant quaternion matrix of size $mp\times np$, at the cost of $O(mnp\log p)$ via the FFT. As a result, we
propose a fast algorithm to calculate the T-product between $m\times n\times p$ and $n\times s\times p$ third-order quaternion tensors via FFTs, at the cost of $O(mnsp)$, which is almost $1/p$ of the computational magnitude of computing T-product by its definition. Numerical calculations verify the correctness of the complexity analysis.

\vspace{3mm}

\noindent {\bf Key words:}\hspace{2mm} block circulant quaternion matrix, block diagonalization, fast Fourier transform, quaternion tensor, tensor-tensor product. \vspace{3mm}

\noindent {\bf Mathematics Subject Classifications (2010):}\hspace{2mm} 15A69, 15B05, 15B33. \vspace{3mm}

\end{abstract}

\section{Introduction}
\label{intro}
It is well-known that block circulant real matrices has been widely applied in various fields, such as code theory \cite{TSSFC-2004}, vibration analysis \cite{OSSPP-2014}, structural calculation \cite{GV-2013,RE-2010,H-1992,CN-1996,CN-2002}, tensor decomposition \cite{D-2016,KM-2011,KMP-2008,ZSK-2018} and image processing \cite{YMYZ-2010,ZEAHK-2014,KHAN-2021,HK-2010}, etc. The wide application of block circulant real matrices greatly thanks to the important property that a block circulant real matrix of size $mp\times np$ can be fast block diagonalization by the discrete Fourier matrix at a cost of $O(mnp\log p)$ via the fast Fourier transform (FFT) \cite{GV-2013,RE-2010}, which brings great convenience for numerical calculations of block circulant real matrices.

One of important applications of block circulant real matrices is the fast calculation of the tensor-tensor product (T-product) between third-order real tensors and T-product based methods for the image and grey video processing. The T-product, proposed by Kilmer and Martin \cite{KM-2011}, is a closed multiplication operation between
third-order tensors, which provides familiar tools of linear algebra for third-order real tensors through a bridge between third-order real tensors and block circulant real matrices \cite{KBHH-2013}. Moreover, with the help of fast block diagonalization of block circulant real matrices via FFTs, the T-product between third-order real tensors is further transformed into block diagonal matrix multiplication in the discrete Fourier domain, which makes tensor factorizations based on the T-product for third-order real tensors can be fast computed and promotes the remarkable success of numerous real tensor methods based on the T-product in the color image and gray video processing \cite{KMP-2008,B-2010,HK-2010,KM-2011,KBHH-2013,ZEAHK-2014,ZA-2017,SNZ-2020,ZBN-2020,Lu-2019,Lu-2016}. However, for real-life problems which are described by fourth-order real tensor models, such as the color video processing, the T-product based method for third-order real tensors is not suitable.

Recently, the quaternion has been found as a powerful tool in the color image and color video processing \cite{SZWGZY-2019,GYWW-2015,HD-2019,JNS-2019,MKL-2020,CXZ-2020,JNS-2019-1,QLWZ-2021,QMZ-2022}. Proposed by Hamilton in 1843 \cite{H-1866}, a quaternion $q$ is often given by:
$
q=a+b{\mathbf{i}}+c{\mathbf{j}}+d{\mathbf{k}},
$ where $a,b,c,d\in \mathbb{R}$, ${\mathbf{i}}$, ${\mathbf{j}}$ and ${\mathbf{k}}$ are three imaginary units which satisfy:
$$\left\{\begin{array}{l}
{\mathbf{i}}^2={\mathbf{j}}^2={\mathbf{k}}^2={\mathbf{i}}{\mathbf{j}}{\mathbf{k}}=-1,\\
{\mathbf{i}}{\mathbf{j}}=-{\mathbf{j}}{\mathbf{i}}={\mathbf{k}},{\mathbf{j}}{\mathbf{k}}
=-{\mathbf{k}}{\mathbf{j}}={\mathbf{i}},
{\mathbf{k}}{\mathbf{i}}=-{\mathbf{i}}{\mathbf{k}}={\mathbf{j}}.
\end{array}\right.$$
It is shown in literatures that the quaternion based expression approach treats three color channel pixels of color images and color videos holistically as a vector field, via encoding the red, green and blue channel pixel values on the three imaginary parts of quaternion matrices or quaternion tensors. As a result, the color structure of color images and color videos can be better maintained in quaternion based methods, compared with the real tensor based methods \cite{MKL-2020,JNS-2019,QLWZ-2021}. Besides, it should be noted that by the quaternion based expression approach, a color image or color video can be described by tensor data which has a lower dimension than the corresponding real tensor data. For example, a color video is described by a fourth-order real tensor, while a color video is described by a third-order quaternion tensor \cite{MKL-2020,JNS-2019,QLWZ-2021}. This observation, together with the good performances of quaternion based methods and T-product based third-order real tensor methods in the color image processing, generates a natural thought, that is, to establish the T-product based third-order quaternion tensor methods for the color video processing.

Analogue to the T-product based third-order real tensor methods, the achievement of the block diagonalization of block circulant quaternion matrices is a key and inevitable issue in the establishment of T-product based third-order quaternion tensor methods, which aims to find a $p\times p$ unitary matrix $\mathbf{P}$ such that for any block circulant quaternion matrix $\mathbf{A}={\rm bcirc}(\mathbf{A_1},\mathbf{A_2},\ldots, \mathbf{A_p})$ of size $mp\times np$,
\begin{equation}\label{equ:1-1}
(\mathbf{P}\otimes \mathbf{I_m})\; \mathbf{A}\;(\mathbf{P}\otimes \mathbf{I_n})^* = {\rm Diag}(\mathbf{\hat{A}_1},\mathbf{\hat{A}_2},\ldots,\mathbf{\hat{A}_p}),
\end{equation}
where ${\rm Diag}(\mathbf{\hat{A}_1},\mathbf{\hat{A}_2},\ldots,\mathbf{\hat{A}_p})$ means the block diagonal matrix with each $\mathbf{\hat{A}_i}$ being its diagonal block entry. 
For a  block circulant real matrix $\mathbf{A}={\rm bcirc}(\mathbf{A_1},\mathbf{A_2},\ldots, \mathbf{A_p})$ with each $\mathbf{A_i}\in \mathbb{R}^{m\times n}$, it has been shown that $\mathbf{A}$ can be block diagonalized by the unitary complex matrices $\mathbf{F_p}\otimes \mathbf{I_m}$ and $\mathbf{F_p}\otimes \mathbf{I_n}$, that is the unitary matrix $\mathbf{P}$ in \eqref{equ:1-1} can be chosen as the DFT matrix $\mathbf{F_p}$ when $\mathbf{A}$ is real. In \cite{QMZ-2022}, it was shown that a type of special block structure quaternion matrices of size $mp\times np$ can block diagonalized by $\mathbf{F_p}\otimes \mathbf{I_m}$ and $\mathbf{F_p}\otimes \mathbf{I_n}$. But for block circulant quaternion matrices, it is not the case. Take the $4\times 4$ circulant quaternion matrix $\mathbf{A}=\operatorname{bcirc}(\mathbf{j},2\mathbf{j},3\mathbf{j},4\mathbf{j})$ for instance. It is easy to deduce
\begin{equation}\label{eq:F4}
\mathbf{F_4}{\mathbf{A}}\mathbf{F_4^*}\\
=\left(
                                             \begin{array}{cccc}
                                               40  &               &        &    \\
                                                   &               &        & -8-8\mathbf{k}\\
                                                   &               & -8     &    \\
                                                   & -8+8\mathbf{k}&        &    \\
                                             \end{array}
                                           \right),
\end{equation}
which implies that unlike the real case, the unitary complex matrices $\mathbf{F_p}\otimes \mathbf{I_m }$ and $\mathbf{F_p}\otimes \mathbf{I_n}$ cannot transform any block circulant quaternion matrices of size $mp\times np$ into block diagonal form. To our best knowledge, the block diagonalization of block circulant quaternion matrices has not been solved ever before, due to that the quaternion field is a noncommutative division algebra \cite{Z-1997}.

Motivated by the above, we aim to find unitary matrix $\mathbf{P}$ such that \eqref{equ:1-1} holds. More importantly, inspired by the fast way via FFTs to calculte the result block diagonal matrix of the $mp\times np$ block circulant real matrix achieved with the help of the DFT matrix $\mathbf{F_p}$, we intend to find such unitary matrices that enables us to inherit the calculation convenience of FFTs to quickly obtain the result block diagonal matrices of block circulant quaternion matrices. The main contributions of this work are as follows:
\begin{itemize}
 \item We show that for any $p\times p$ circulant quaternion matrix $\mathbf{A}$, there exist $\mathbf{P},\mathbf{Q}\in \mathbb{C}^{p\times p}$ such that \begin{equation}\label{eq:Fp}
     (\mathbf{F_p}\mathbf{p}){\mathbf{A}}(\mathbf{F_p}\mathbf{p})^*=
     (\mathbf{F_p}{\mathbf{P}}\mathbf{F^*_p})+(\mathbf{F_p}{\mathbf{Q}}\mathbf{F_p})\mathbf{j},
     \end{equation}
     where $\mathbf{p}=1,\mathbf{j}$ or $(1+\mathbf{j})/\sqrt{2}$, since $1c = c1$ but $\mathbf{j}c = \bar{c}\mathbf{j}$ for any complex number $c$. This means that $\mathbf{F_p}\mathbf{p}$ with $\mathbf{p}=1,\mathbf{j}$ or $(1+\mathbf{j})/\sqrt{2}$ cannot be diagonalization matrix of $\mathbf{A}$, since the part of $(\mathbf{F_p}{\mathbf{Q}}\mathbf{F_p})$ is not diagonal in general as shown in \eqref{eq:F4}. Then we prove that to find the unitary quaternion diagonalization matrix of circulant quaternion matrices is equivalent to solve eight diagonalization problem of real circulant quaternion matrices, which is hard to achieve.
 \item We provide several special cases of the multiplication between special octonions and quaternions to avoid the ``$\mathbf{F_p}{\mathbf{Q}}\mathbf{F_p}$'' part in \eqref{eq:Fp}, observing that $\mathbf{l}q=\bar{q}\mathbf{l}$ and $(p\mathbf{l})(q\mathbf{l})$ is a quaternion for any $p,q$ being quaternions. With the help of that, we show that for any $p\times p$ circulant quaternion matrix with the unique Cayley\^{a}Dickson (CD) form $\mathbf{A}=\mathbf{A_1}+\mathbf{A_2}\mathbf{j}$ where $\mathbf{A_1},\mathbf{A_2}\in \mathbb{C}^{p\times p}$, $$(\mathbf{F_p}\mathbf{p}){\mathbf{A}}(\mathbf{F_p}\mathbf{p})^*=
     (\mathbf{F_p}\overline{{\mathbf{A_1}}}\mathbf{F^*_p})-
     (\mathbf{F_p}{\mathbf{A_2}}\mathbf{F^*_p})\mathbf{j},$$ where $\mathbf{p}=\mathbf{l},\mathbf{jl}$ or $(\mathbf{l}+\mathbf{jl})/\sqrt{2}$. That means, $\mathbf{F_p}\mathbf{p}$ with $\mathbf{p}=\mathbf{l},\mathbf{jl}$ or $(\mathbf{l}+\mathbf{jl})/\sqrt{2}$ is a diagonalization matrix of $p\times p$ circulant quaternion matrix, thus we obtain a fast way to compute the result diagonal matrix of $p\times p$ circulant quaternion matrices via FFTs at the cost of $O(p\log p)$. Furthermore, we prove that each $mp\times np$ block circulant quaternion matrix can be transformed into block diagonal quaternion matrix by $\mathbf{F_p}\mathbf{p}\otimes \mathbf{I_m }$ and $\mathbf{F_p}\mathbf{p}\otimes \mathbf{I_n}$ via FFTs at the cost of $O(mnp\log p)$.
 \item We establish a relationship between T-product of third-order quaternion tensors and the block diagonalization of block circulant quaternion matrices, and provide a way to compute the T-product of two third-order quaternion tensors of size $m\times n\times p$ and $n\times s\times p$ via FFTs at the cost of $O(mnsp)$. For third-order quaternion tensors $\mathcal{A}$ of size $m\times n\times p$ and $\mathcal{B}$ of size $n\times s\times p$, the T-product $\mathcal{C}=\mathcal{A}\ast\mathcal{B}$ can be obtained by the following procedure:
     \begin{itemize}
     	 \item First, compute $\hat{\mathcal{A}}=\operatorname{fft}(\overline{\mathcal{A}},[\;],3)$ and $\hat{\mathcal{B}}=\operatorname{fft}(\overline{\mathcal{B}},[\;],3)$, denoted by the CD forms $\hat{\mathcal{A}}=\hat{\mathcal{A}_1}+\hat{\mathcal{A}_2}\mathbf{j}$ and $\hat{\mathcal{B}}=\hat{\mathcal{B}_1}+\hat{\mathcal{B}_2}\mathbf{j}$ ; \item Second, let $\hat{\mathcal{A}}(:,:,p+1):=\hat{\mathcal{A}}(:,:,1)$ and  $\hat{\mathcal{B}}(:,:,p+1):=\hat{\mathcal{B}}(:,:,1)$, compute
     \begin{eqnarray*}
     {\hat{\mathcal{C}}}_1(:,:,i)&=& {\hat{\mathcal{A}}_1}(:,:,p+2-i)\hat{\mathcal{B}}_1(:,:,p+2-i)
     -\overline{{\hat{\mathcal{A}}_2}}(:,:,i)\hat{\mathcal{B}}_2(:,:,p+2-i), \\
     {\hat{\mathcal{C}}_2}(:,:,i)&=&\overline{\hat{\mathcal{A}}_1}(:,:,p+2-i)
     \hat{\mathcal{B}}_2(:,:,i)+\hat{\mathcal{A}}_2(:,:,i)\hat{\mathcal{B}}_1(:,:,i); \end{eqnarray*}
\item Last, compute $\mathcal{C}=\operatorname{ifft}(\overline{\hat{\mathcal{C}}},[\;],3)$ where $\hat{\mathcal{C}}=\hat{\mathcal{C}_1}+\hat{\mathcal{C}_2}\mathbf{j}$.
     \end{itemize}
    Numerical experiments demonstrate that the proposed algorithm speeds up $p$ times of the calculation of the T-product of two third-order quaternion tensors with $p$ frontal slices.
 \item The achieved block diagonalization of block circulant quaternion opens up a way to extend the greatly successful T-product based real tensor methods to the quaternion case, meanwhile this paper provides a new vein for generalizing more matrix theory in the complex domain to the quaternion domain through the participation of octonion matrices.
\end{itemize}

The organization of this paper is divided into six sections. In Section 2, we introduce preliminary knowledge, including some notations, octonions, the block diagonalization of block circulant complex matrices, and T-product of third-order tensors. In Section 3, we study the diagonalization of circulant quaternion matrices in the quaternion domain. We show that $\mathbf{F_p}\mathbf{p}$ with $\mathbf{p}=1,\mathbf{j}$ or $(1+\mathbf{j})/\sqrt{2}$ is not diagonalization matrix of $p\times p$ circulant quaternion matrices and give sufficient and necessary conditions for a quaternion matrix being a quaternion diagonalization matrix of circulant quaternion matrices. In Section 4, we study the diagonalization of circulant quaternion matrices in the octonion domain. We provide several special cases of the non-commutative and non-associative octonion multiplication. Then, we prove that $\mathbf{F_p}\mathbf{p}$ with $\mathbf{p}=\mathbf{l},\mathbf{jl}$ or $(\mathbf{l}+\mathbf{jl})/\sqrt{2}$ is unitary and is a diagonalization matrix of $p\times p$ circulant quaternion matrix. Finally, we achieve the block diagonalization of $mp\times np$ block circulant quaternion matrices with octonion matrices $\mathbf{F_p}\mathbf{p}\otimes \mathbf{I_m } $ and $\mathbf{F_p}\mathbf{p}\otimes \mathbf{I_n}$, where $\mathbf{p}=\mathbf{l},\mathbf{jl}$ or $(\mathbf{l}+\mathbf{jl})/\sqrt{2}$.
In Section 5, we give an algorithm to fast compute the T-product of quaternion tensors by applying the fast block diagonalization way of block circulant quaternion matrices via FFTs. Last, the concluding remarks are made in Section 6. 
\section{Preliminary}
\label{Sect. 1}

\subsection{Notation}
\label{Sect. 1.1}
Throughout this paper, we use small letters $a,b,\ldots$ for scalars, small bold letters $\mathbf{{a}}, \mathbf{{b}}, \ldots$ for vectors, bold capital letters $\mathbf{A},\mathbf{B},\ldots$ for matrices, and calligraphic letters $\mathcal{A}, \mathcal{B}, \ldots$ for tensors. For any positive integer $n$, denote $[n] := \{1,2,\ldots , n\}$. Denote the set of all real numbers, complex numbers, quaternions and octonions as $\mathbb{R}$, $\mathbb{C}$, $\mathbb{Q}$ and $\mathbb{O}$, respectively. In the following, we use $\mathbb{F}$ to represent $\mathbb{R}$, $\mathbb{C}$, $\mathbb{Q}$ or $\mathbb{O}$, respectively, unless otherwise specified.
Denote $\mathbb{F}^{m\times n}$ by the collection of all $m\times n$ matrices with entries in $\mathbb{F}$. When $n=1$, we simplify $\mathbb{F}^{m\times n}$ as $\mathbb{F}^{m}$. In addition, we use $\mathbf{I_m}$ to represent the $m\times m$ identity matrix.
For any ${\mathbf{{X}}}:=(x_{ij})\in \mathbb{F}^{m\times n}$, let $\overline{{{\mathbf{X}}}}:=(\overline{x_{ij}})$ and ${{\mathbf{X}}}^{*}:=(\overline{x_{ij}})^\top$ represent the conjugate and conjugate transpose of the matrix ${{X}}$, respectively. For any $\mathbf{x}:=(x_r)\in \mathbb{F}^n$, ${\rm Diag}(\mathbf{x})$ means the diagonal matrix whose diagonal entries are $x_r$, respectively. For any $\mathbf{{X_r}}\in \mathbb{F}^{m\times n}$ with $r\in[p]$, ${\rm Diag}(\mathbf{X_1}, \mathbf{X_2},\ldots, \mathbf{X_p}) :={\rm Diag}({\mathbf{X_r}}:r\in[p])$ means the block diagonal matrix whose diagonal block entries are ${\mathbf{X_r}}$, respectively. For any matrices $\mathbf{A}$ and $\mathbf{B}$ of appropriate sizes, $\mathbf{A}\cdot\mathbf{B}$ denotes the product between $\mathbf{A}$ and $\mathbf{B}$, and sometimes $\cdot$ is omitted. For any third-order tensors $\mathcal{A}$ and $\mathcal{B}$ of appropriate sizes, $\mathcal{A}\ast\mathcal{B}$ denotes the T-product between $\mathcal{A}$ and $\mathcal{B}$.


\subsection{Octonions and quaternions}
\label{Sect. 1.2}
Octonions were first discovered by John Graves, and later independently introduced by Cayley in 1845 \cite{C-1845}. An octonion $o\in \mathbb{O}$ is often given by:
$o=a+b{\mathbf{i}}+c{\mathbf{j}}+d{\mathbf{k}}+e{\mathbf{l}}
+f{\mathbf{il}}+g{\mathbf{jl}}+h{\mathbf{kl}},
$
where $a,b,c,d,e,f,g,h\in \mathbb{R}$, ${\mathbf{i}}$, ${\mathbf{j}}$, ${\mathbf{k}}$, ${\mathbf{l}}$, ${\mathbf{il}}$, ${\mathbf{jl}}$ and ${\mathbf{kl}}$ are square roots of $-1$, which satisfy the following multiplication stable:

\begin{table}[H]
  \centering
\begin{tabular}{c|cccccccc}

  $\mathbb{O}$ & ${1}$ & ${\mathbf{i}}$ & ${\mathbf{j}}$ & ${\mathbf{k}}$ & ${\mathbf{l}}$ & ${\mathbf{il}}$ & ${\mathbf{jl}}$ & ${\mathbf{kl}}$ \\  \hline
  ${1}$ & ${1}$ & ${\mathbf{i}}$ & ${\mathbf{j}}$ & ${\mathbf{k}}$ & ${\mathbf{l}}$ & ${\mathbf{il}}$ & ${\mathbf{jl}}$ & ${\mathbf{kl}}$  \\
  ${\mathbf{i}}$ & ${\mathbf{i}}$ & -$1$ & ${\mathbf{k}}$ & -${\mathbf{j}}$ & ${\mathbf{il}}$ & -${\mathbf{l}}$ & -${\mathbf{kl}}$ & ${\mathbf{jl}}$ \\
  ${\mathbf{j}}$ & ${\mathbf{j}}$ & -${\mathbf{k}}$ & -$1$ & ${\mathbf{i}}$ & ${\mathbf{jl}}$ & ${\mathbf{kl}}$ & -${\mathbf{l}}$ & -${\mathbf{il}}$ \\
  ${\mathbf{k}}$ & ${\mathbf{k}}$ & ${\mathbf{j}}$ & -${\mathbf{i}}$ & -$1$ & ${\mathbf{kl}}$ & -${\mathbf{jl}}$ & ${\mathbf{il}}$ & -${\mathbf{l}}$ \\
  ${\mathbf{l}}$ & ${\mathbf{l}}$ & -${\mathbf{il}}$ & -${\mathbf{jl}}$ & -${\mathbf{kl}}$ & -$1$ & ${\mathbf{i}}$ & ${\mathbf{j}}$ & ${\mathbf{k}}$ \\
  ${\mathbf{il}}$ & ${\mathbf{il}}$ & ${\mathbf{l}}$ & -${\mathbf{kl}}$ & ${\mathbf{jl}}$ & -${\mathbf{i}}$ & -$1$ & -${\mathbf{k}}$ & ${\mathbf{j}}$ \\
  ${\mathbf{jl}}$ & ${\mathbf{jl}}$ & ${\mathbf{kl}}$ & ${\mathbf{l}}$ & -${\mathbf{il}}$ & -${\mathbf{j}}$ & ${\mathbf{k}}$ & -$1$ & -${\mathbf{i}}$ \\
  ${\mathbf{kl}}$ & ${\mathbf{kl}}$ & -${\mathbf{jl}}$ & ${\mathbf{il}}$ & ${\mathbf{l}}$ & -${\mathbf{k}}$ & -${\mathbf{j}}$ & ${\mathbf{i}}$ &  -$1$\\
\end{tabular}
  \caption{Octonion Multiplication Rule Table}\label{multiplication-octonion}
\end{table}

For any given two octonions $o_1=a_1+b_1{\mathbf{i}}+c_1{\mathbf{j}}+d_1{\mathbf{k}}+e_1{\mathbf{l}}+f_1{\mathbf{il}}+g_1{\mathbf{jl}}+h_1{\mathbf{kl}},
o_2=a_2+b_2{\mathbf{i}}+c_2{\mathbf{j}}+d_2{\mathbf{k}}+e_2{\mathbf{l}}+f_2{\mathbf{il}}+g_2{\mathbf{jl}}+h_2{\mathbf{kl}}\in \mathbb{O}$, the sum is defined as:
$$\begin{array}{rl}
o_1+o_2:=&(a_1+a_2)+(b_1+b_2){\mathbf{i}}+(c_1+c_2){\mathbf{j}}+(d_1+d_2){\mathbf{k}}\\
&+(e_1+e_2){\mathbf{l}}+(f_1+f_2){\mathbf{il}}+(g_1+g_2){\mathbf{jl}}+(h_1+h_2){\mathbf{kl}},
\end{array}$$
and the multiplication of $o_1$ and $o_2$ is defined as:
$$\begin{array}{rccll}
o_1o_2&:=& &(a_1a_2-b_1b_2-c_1c_2-d_1d_2-e_1e_2-f_1f_2-g_1g_2-h_1h_2)&\\
  & &+&(a_1b_2+b_1a_2+c_1d_2-d_1c_2+e_1f_2-f_1e_2+h_1g_2-g_1h_2)&{\mathbf{i}}\\
  & &+&(a_1c_2+c_1a_2+d_1b_2-b_1d_2+e_1g_2-g_1e_2+f_1h_2-h_1f_2)&{\mathbf{j}}\\
  & &+&(a_1d_2+d_1a_2+b_1c_2-c_1b_2+e_1h_2-h_1e_2+g_1f_2-f_1g_2)&{\mathbf{k}}\\
  & &+&(a_1e_2+e_1a_2+f_1b_2-b_1f_2+g_1c_2-c_1g_2+h_1d_2-d_1h_2)&{\mathbf{l}}\\
  & &+&(a_1f_2+f_1a_2+b_1e_2-e_1b_2+h_1c_2-c_1h_2+d_1g_2-g_1d_2)&{\mathbf{il}}\\
  & &+&(a_1g_2+g_1a_2+b_1h_2-h_1b_2+c_1e_2-e_1c_2+f_1d_2-d_1f_2)&{\mathbf{jl}}\\
  & &+&(a_1h_2+h_1a_2+g_1b_2-b_1g_2+c_1f_2-f_1c_2+d_1e_2-e_1d_2)&{\mathbf{kl}},
\end{array}$$
by Table \ref{multiplication-octonion}. For any $o\in \mathbb{O}$, its conjugate $\bar{o}$ and modulus $|o|$ are, respectively, defined as:
\begin{eqnarray}
  \label{conj}
\bar{o} &:=& a-b{\mathbf{i}}-c{\mathbf{j}}-d{\mathbf{k}}-e{\mathbf{l}}-f{\mathbf{il}} -g{\mathbf{jl}}-h{\mathbf{kl}},  \\
\nonumber  |o| &:=& \sqrt{\bar{o}o} = \sqrt{a^2+b^2+c^2+d^2+e^2+f^2+g^2+h^2}.
\end{eqnarray}

When $e,f,g,h$ are zeros, the octonion $o=a+b{\mathbf{i}}+c{\mathbf{j}}+d{\mathbf{k}}+e{\mathbf{l}}
+f{\mathbf{il}}+g{\mathbf{jl}}+h{\mathbf{kl}}$ becomes a quaternion. Obviously, each octonion $o=a+b{\mathbf{i}}+c{\mathbf{j}}+d{\mathbf{k}}+e{\mathbf{l}}+f{\mathbf{il}}+g{\mathbf{jl}}+h{\mathbf{kl}}$ has a uniquely quaternionic form as $o=q_1+q_2{\mathbf{l}}$, where $q_1=a+b{\mathbf{i}}+c{\mathbf{j}}+d{\mathbf{k}}\in \mathbb{Q}$ and $q_2=e+f{\mathbf{i}}+g{\mathbf{j}}+h{\mathbf{k}}\in \mathbb{Q}$. Similarly, each quaternion $q=a+b{\mathbf{i}}+c{\mathbf{j}}+d{\mathbf{k}}$ can be uniquely expressed by the CD form, i.e., the complex form $q=p_1+p_2{\mathbf{j}}$, where $p_1=a+b{\mathbf{i}}\in \mathbb{C}$ and $p_2=c+d{\mathbf{i}}\in \mathbb{C}$. It is easy to deduce that for quaternions, the associative law holds, but the commutative law does not hold, and for octonions, both the associative law and commutative law do not hold, i.e.,
\begin{eqnarray*}
  && q_1(q_2q_3)= (q_1q_2) q_3,\ q_1q_2\not= q_2q_1,\ {\rm for\ general}\ q_1, q_2, q_3\in \mathbb{Q},  \\
  && o_1(o_2o_3)\not= (o_1o_2) o_3,\ o_1o_2\not= o_2o_1,\ {\rm for\ general}\  o_1,o_2,o_3\in \mathbb{O}.
\end{eqnarray*}

Besides, an $m\times n$ matrix $\mathbf{A}$ is said a quaternion (or octonion, respectively) matrix if its entries are quaternions (or octonion, respectively), i.e., $\mathbf{A}_{ij}\in \mathbb{Q}$
(or $\mathbf{A}_{ij}\in\mathbb{O}$, respectively) for all $i\in [m],j\in [n]$.
A square $n\times n$ quaternion (or octonion, respectively) matrix $\mathbf{A}$ is said to be inversible if there exists some quaternion (or octonion, respectively) matrix $\mathbf{B}$ such that $\mathbf{AB}=\mathbf{BA}=\mathbf{I_n}$. Then the matrix $\mathbf{B}$ is called an inverse of $\mathbf{A}$, denoted by $\mathbf{A}^{-1}$. In addition, a quaternion (or octonion, respectively) matrix $\mathbf{A}$ is called unitary if $\mathbf{A A^*}=\mathbf{A^*A}=\mathbf{I_n}.$

%
%

\subsection{Block diagonalization of block circulant complex matrices}
\label{Sect.2}
In this subsection, we recall the block diagonalization of block circulant complex matrices. Let $\mathbb{F}= \mathbb{R}, \mathbb{C}, \mathbb{Q}$ or $\mathbb{O}$ be a number field. A matrix ${\mathbf{A}}\in \mathbb{F}^{mp\times np}$ is called a block circulant matrix if it has the following form \cite{D-1994}:
\begin{equation}\label{eq:circulant}
\mathbf{{A}}=\left(\begin{array}{ccccc}
\mathbf{{A_1}}& \mathbf{{A_p}} & \cdots & \mathbf{{A_3}}& \mathbf{{A_2}}\\
\mathbf{{A_2}}& \mathbf{{A_1}} & \cdots & \mathbf{{A_4}}& \mathbf{{A_3}}\\
\vdots & \vdots & \ddots  & \vdots & \vdots\\
\mathbf{{A_p}} & \mathbf{{A_{p-1}}} & \cdots & \mathbf{{A_2}}& \mathbf{{A_1}}
\end{array}\right)
\end{equation}
where $\mathbf{{A_r}}\in \mathbb{F}^{m\times n}$ for all $r\in [p]$.
For simplicity, $A$ is often written as $\mathbf{{A}}=\operatorname{bcirc}(\mathbf{{A_1}},\mathbf{{A_2}},\ldots, \mathbf{{A_p}}).$
In particular, $\mathbf{{A}}$ as in \eqref{eq:circulant} is called a circulant matrix if $m=n=1$ and denoted by $\mathbf{{A}}=\operatorname{circ}(a_1,a_2,\ldots, a_p).$

In the theory of diagonalization of circulant matrices, the diagonalization matrix plays an important role, which is a unitary matrix defined as follow.
\begin{definition}(Diagonalization matrix)
Let $\mathbf{F}\in \mathbb{F}^{p\times p}$ be a unitary matrix. If $\mathbf{FAF^*}$ is a diagonal matrix for any $p\times p$ circulant matrix $\mathbf{A}$, then $\mathbf{F}$ is called a diagonalization matrix of $p\times p$ circulant matrices.
\end{definition}

Denote $\mathbf{F_{p}}$ as the $p\times p$ normalized discrete Fourier transform (DFT) matrix \cite{GV-2013}, i.e.,
\begin{eqnarray}\label{eq:DFT}
\mathbf{F_{p}}=\frac{1}{\sqrt{p}}\left(
\begin{array}{ccccc}
  1 & 1 & 1  & \ldots & 1 \\
  1 & \omega & \omega^2 & \ldots & \omega^{p-1} \\
  1 & \omega^2 & \omega^4  & \ldots & \omega^{2(p-1)} \\
  \vdots & \vdots & \vdots  & \ddots & \vdots \\
  1 & \omega^{p-1} & \omega^{2(p-1)} & \ldots & \omega^{(p-1)(p-1)}
\end{array}
\right)
\end{eqnarray}
where $\omega=e^{-\frac{2\pi \mathbf{i}}{p}}\in \mathbb{C}$. It is deduced that
\begin{equation}\label{FpFp}
\mathbf{F_p}\mathbf{F_p}=\mathbf{F_p^*}\mathbf{F_p^*}=\mathbf{P_p},
\end{equation}
where $\mathbf{P_p}$ is a permutation matrix defined by
\begin{equation}\label{eq:Pp}
\left\{\begin{array}{ll}
P_{11}=P_{ij}=1, & {\rm if } \quad i+j=p+2,\quad i,j=2,3,\ldots,p;\\
P_{ij}=0, & {\rm otherwise }.\\
\end{array}\right.
\end{equation}

Recall that the Kronecker product $\mathbf{A}\otimes \mathbf{B}$, of a matrix $\mathbf{A}\in \mathbb{F}^{m\times n}$ and a matrix $\mathbf{B}\in \mathbb{F}^{s\times t}$ is an $ms\times nt$ block matrix defined by
\begin{equation}\label{kronecker}
\mathbf{A}\otimes \mathbf{B}=\left(\begin{array}{cccc}
a_{11} \mathbf{B}& a_{11}\mathbf{ B} & \ldots & a_{1n} \mathbf{B }\\
a_{21}\mathbf{ B} & a_{22} \mathbf{B} & \ldots & a_{2n} \mathbf{B} \\
\vdots & \vdots &  \ddots & \vdots \\
a_{m1} \mathbf{B} & a_{m2} \mathbf{B} & \cdots & a_{mn} \mathbf{B}
\end{array}
\right).
\end{equation}
For positive integers $p$ and $m$, it is deduced directly that
\begin{equation}\label{eq:invtranFpIm}
 (\mathbf{F_p}\otimes \mathbf{I_m})^* = \mathbf{F}^*_p\otimes \mathbf{I_m}.
\end{equation}


It is known that any circulant complex matrix ${\mathbf{A}}=\operatorname{circ}(a_1,a_2,\ldots, a_p)\in \mathbb{C}^{p\times p}$ can be diagonalized by the $p\times p$ DFT matrix $\mathbf{F_{p}}$ at the cost of $O(p\log p)$ flops, see ref. \cite{GV-2013}, i.e.,
\begin{equation}\label{eq:diagonal-C}
\mathbf{F_{p}}\mathbf{{A}}\mathbf{F_{p}^*}=\mathbf{F_{p}}[\operatorname{circ}(a_1,a_2,\ldots, a_p)]\mathbf{F_{p}^*}={\rm Diag}(d_1,d_2,\ldots,d_p)
\in \mathbb{C}^{p\times p}.
\end{equation}
That is to say, $\mathbf{F_p}$ is a diagonalization matrix of $p\times p$ circulant complex matrices. Besides, the diagonal entries $d_1,d_2,\ldots,d_p$ of the diagonal matrix in \eqref{eq:diagonal-C} and the generators $a_1,a_2,\ldots, a_p$ of the circulant matrix $\mathbf{A}$ in \eqref{eq:diagonal-C} satisfy the following relations:
$$(d_1,d_2,\ldots, d_p)^\top=\sqrt{p}\mathbf{F_{p}}(a_1,a_2,\ldots, a_p)^\top.$$ Thus, the vector $d=(d_1,d_2,\ldots, d_p)^\top$ can be computed by applying the fast fourier transform (FFT for short) to the vector $a=(a_1,a_2,\ldots, a_p)^\top$, i.e.,
$d=\operatorname{fft}(a).$

Furthermore, any block circulant complex matrix $\mathbf{{A}}\in \mathbb{C}^{mp\times np}$ defined as in (\ref{eq:circulant}) can be block diagonalized by $\mathbf{F_p}\otimes \mathbf{I_m }\in \mathbb{C}^{mp\times mp}$
and $\mathbf{F_p}\otimes \mathbf{I_n} \in \mathbb{C}^{np\times np}$
as
$$(\mathbf{F_p}\otimes \mathbf{I_m}) \mathbf{{A}}(\mathbf{F^*_p}\otimes \mathbf{I_n}) ={\rm Diag}(\mathbf{\hat{A}_1},\mathbf{\hat{A}_2},\ldots,\mathbf{\hat{A}_p})\in \mathbb{C}^{mp\times np}, $$
where $\mathbf{\hat{A}_r}\in \mathbb{C}^{m\times n}$ for any $r\in [p]$.

Then, the following questions naturally arise.

\noindent{\bf Question 1}: Can the unitary complex matrices $\mathbf{F_p}\otimes \mathbf{I_m}$
and $\mathbf{F_p}\otimes \mathbf{I_n}$ block diagonalize an arbitrary block circulant quaternion matrix $\mathbf{{A}}\in \mathbb{Q}^{mp\times np}$? And can the unitary quaternion matrices $(\mathbf{F_p}\mathbf{p})\otimes \mathbf{I}_m$
and $(\mathbf{F_p}\mathbf{p})\otimes \mathbf{I_n}$ with $\mathbf{p}=\mathbf{j}$ or $(1+\mathbf{j})/\sqrt{2}$ block diagonalize an arbitrary block circulant quaternion matrix ${A}\in \mathbb{Q}^{mp\times np}$?

\noindent{\bf Question 2}: If the answer to {\bf Question 1} is negative, then what kind of unitary quaternion matrices $\mathbf{P}\otimes \mathbf{I_m}\in \mathbb{Q}^{mp\times mp}$ and $\mathbf{P}\otimes \mathbf{I_n}\in \mathbb{Q}^{np\times np}$ can diagonalize an arbitrary block circulant quaternion matrix $\mathbf{{A}}\in \mathbb{Q}^{mp\times np}$ ? i.e.,
$
(\mathbf{P}\otimes \mathbf{I_m})\mathbf{A}(\mathbf{P}\otimes \mathbf{I_n})^* = {\rm Diag}(\mathbf{\hat{A}_1},\mathbf{\hat{A}_2},\ldots,\mathbf{\hat{A}_p}).
$

These questions above are answered in Section \ref{Sect.DQCM}. We give negative answers to {\bf Question 1} and prove sufficient and necessary conditions for a unitary quaternion matrix being a quaternion diagonalization matrix of circulant quaternion matrices.  
\subsection{T-product of third-order tensors}
\label{Sub:T-product}
Here, we recall some preliminaries for the definition of T-product between third-order tensors and show that how T-product between third-order real tensors can be fast calculated with the help of the block diagonalization way of block circulant real matrices.

A third-order tensor $\mathcal{A}\in \mathbb{F}^{m\times n\times p}$ can be treated as a stack of frontal slices. The $r$-th frontal slice is denoted as $\mathbf{{A}^{(r)}}\in \mathbb{F}^{m\times n}$ for any $r\in [p]$, i.e., $\mathbf{{A}^{(r)}} = \mathcal{A}(:,:,r).$ A third-order tensor $\mathcal{A}$ is closely related to a block circulant matrix and a block column vector by the ``${bcirc}$'' and ``${unfold}$'' operators, respectively, which are defined as:
\begin{equation}\label{bcirc-unfold}\begin{array}{rcl}
bcirc(\mathcal{A}):=
\left(\begin{array}{cccc}
\mathbf{{A}^{(1)}} & \mathbf{{A}^{(p)}} & \ldots & \mathbf{{A}^{(2)}} \\
\mathbf{{A}^{(2)}} & \mathbf{{A}^{(1)}} & \ldots & \mathbf{{A}^{(3)}} \\
\vdots & \vdots &  \ddots & \vdots \\
\mathbf{{A}^{(p)}} & \mathbf{{A}^{(p-1)}} & \cdots & \mathbf{{A}^{(1)}}
\end{array}
\right)&
,& unfold(\mathcal{A}):=
\left(\begin{array}{c}
\mathbf{ {A}^{(1)}}\\
 \mathbf{{A}^{(2)}}\\
 \vdots\\
 \mathbf{{A}^{(p)}}
\end{array}
\right),
\end{array}
\end{equation}
and their inverse ``$bcirc^{-1}$" and ``$fold$" are defined as:
$$ bcirc^{-1}(bcirc(\mathcal{A})):=\mathcal{A},\ fold(unfold(\mathcal{A})):=\mathcal{A}.$$

When $\mathbb{F}=\mathbb{R}$, an important multiplication between third-order tensors whose entries belong to $\mathbb{F}$ was shown in \cite{KM-2011,KBHH-2013}, named {\it T-product}. Here, we apply the same definition of multiplication between third-order tensors for the case of $\mathbb{F}=\mathbb{C},\ \mathbb{Q}$ or $ \mathbb{O}$.

\begin{definition}\label{T-product}
(T-product) Let $\mathcal{A}\in \mathbb{F}^{m\times n\times p}$, $\mathcal{B}\in \mathbb{F}^{n\times s\times p}$. Then the T-product $\mathcal{A}\ast\mathcal{B}$ is defined as:
$$\mathcal{A}\ast\mathcal{B}:=fold(bcirc(\mathcal{A})\cdot unfold(\mathcal{B}))\in \mathbb{F}^{m\times s\times p}.$$
\end{definition}

For the convenience of description, we introduce the following operator:
\begin{equation}\label{eqdef:blockdiag}
{\rm Diag}(\mathcal{A}):=\left(\begin{array}{cccc}
                             \mathbf{{A}^{(1)}} &  &  &  \\
                              & \mathbf{{A}^{(2)}} &  &  \\
                              &  & \ddots &  \\
                              &  &  & \mathbf{{A}^{(p)}}
                           \end{array}\right),
\end{equation}
where $\mathcal{A}$ is a tensor in $\mathbb{F}^{m\times n\times p}$.

As shown in the literatures \cite{KMP-2008,B-2010,KM-2011,KBHH-2013}, for a third-order real tensor, the following results hold.

\begin{proposition} \cite{KM-2011}\label{pro:block-cir}
Let $\mathcal{A}\in \mathbb{F}^{m\times n\times p}$, where $\mathbb{F}=\mathbb{R}$. Then the block circulant matrix $bcirc(\mathcal{A})$ 
can be block diagonalized by $\mathbf{F_p}\otimes \mathbf{I_m }\in \mathbb{C}^{mp\times mp}$ and $\mathbf{F_p}\otimes \mathbf{I_n} \in \mathbb{C}^{np\times np}$, that is
$$({\mathbf{F_p}}\otimes \mathbf{I_m})\ bcirc(\mathcal{A})\ (\mathbf{F^*_p}\otimes \mathbf{I_n})={\rm Diag}(\hat{\mathcal{A}})
$$
and $\hat{\mathcal{A}}$  and $\mathcal{A}$ satisfy:
$$
unfold(\hat{\mathcal{A}})=\sqrt{p}(\mathbf{F_p}\otimes \mathbf{I}_m)unfold(\mathcal{A}).
$$
That means, the tensor $\hat{\mathcal{A}}$ can be computed at the cost of $O(mnp\log p)$, by applying FFT along the third dimension of $\mathcal{A}$, i.e., $\hat{\mathcal{A}}=\operatorname{fft}(\mathcal{A},[\; ],3)$.
\end{proposition}


By the block diagonalization of block circulant complex matrices shown in Proposition \ref{pro:block-cir}, a fast way to calculate the T-product of two third-order real tensors was given as follow \cite{KM-2011}:

\begin{center}
{\begin{tabular}{rl}
  \hline\noalign{\smallskip}
{\bf Algorithm}&{\bf 2.1: Computing the T-product of real tensors via FFTs} \\
\noalign{\smallskip}\hline\noalign{\smallskip}
{\bf \;\;Input:}& $\mathcal{A}\in \mathbb{R}^{m\times n\times p}$ and $\mathcal{B}\in \mathbb{R}^{n\times s \times p}$.\\
{\;\;\bf Output:}& $\mathcal{C}\in \mathbb{R}^{m\times s\times p}$: $\mathcal{C}=\mathcal{A}\ast\mathcal{B}$.\\
{\;\;\bf Step 1:}&Compute $\hat{\mathcal{A}}=\operatorname{fft}(\mathcal{A},[\;],3)$ and $\hat{\mathcal{B}}=\operatorname{fft}(\mathcal{B},[\;],3)$.\\
{\;\;\bf Step 2:}& For $r=1:p,$ compute $r$-th frontal slice of $\hat{\mathcal{C}}$ by\\
&\quad\quad$\hat{\mathcal{C}}(:,:,r)=\hat{\mathcal{A}}(:,:,r)\hat{\mathcal{B}}(:,:,r),$\\
&end for.\\
{\;\;\bf Step 3:}& Compute $\mathcal{C}=\operatorname{ifft}(\hat{\mathcal{C}},[\;],3)$.\\
\noalign{\smallskip}  \hline
\end{tabular}}
\end{center}

Then, the following question arises:

\noindent{\bf Question 3}: When $\mathbb{F}=\mathbb{Q}$, can a similar way as Algorithm 2.1 be given to accelerate the computation of the T-product of third-order quaternion tensors defined by Definition \ref{T-product}?

It is achieved by studies in Section \ref{Sect.3.2} and Section \ref{Sect.4}. In Section \ref{Sect.3.2}, we show how to obtain the result block diagonal quaternion matrix of an arbitrary $mp\times np$ block circulant quaternion matrix via FFTs, which is block diagonalized by unitary octonion matrices $\mathbf{F_p}\mathbf{p}$ with $\mathbf{p}=\mathbf{l},\mathbf{il}$ or $(\mathbf{l}+\mathbf{il})/\sqrt{2}$. Then through establishing the relationship between the T-product of third-order quaternion tensors and the block diagonalization of block circulant quaternion matrices, we provide a similar way as Algorithm 2.1 for the calculation of the T-product of third-order quaternion tensors with $p$ frontal slices in Section \ref{Sect.4}, which speed up almost $p$ times compared with the T-product computation by its definition. 

\section{Diagonalization matrices in the quaternion domain}
\label{Sect.DQCM}
When all entries of a circulant matrix defined by \eqref{eq:circulant} are quaternions, it is called a circulant quaternion matrix. To answer {\bf Question 1} and {\bf Question 2}, in this section we study the diagonalization of circulant quaternion matrices in the quaternion domain, which is to find a unitary quaternion matrix $\mathbf{{F}}\in \mathbb{Q}^{p\times p}$ such that $\mathbf{{F}{A}{F}^{*}}$ is a diagonal matrix for all $\mathbf{{A}}=\operatorname{circ}(a_1,a_2,\ldots, a_p)\in \mathbb{Q}^{p\times p}$, and the unitary quaternion matrix $\mathbf{{F}}$ is called a diagonalization matrix of ${p\times p}$ circulant quaternion matrices.

First, we answer {\bf Question 1} by the following example, which shows that $\mathbf{F_p}\mathbf{p}$ with $\mathbf{p}=1,\mathbf{j}$ or $(1+\mathbf{j})/\sqrt{2}$ cannot diagonalize an arbitrary circulant quaternion matrix ${A}\in \mathbb{Q}^{p\times p}$. Therefore, $(\mathbf{F_p}\mathbf{p})\otimes \mathbf{I}_m$
and $(\mathbf{F_p}\mathbf{p})\otimes \mathbf{I_n}$ with $\mathbf{p}=\mathbf{j}$ or $(1+\mathbf{j})/\sqrt{2}$ cannot block diagonalize an arbitrary block circulant quaternion matrix $\mathbf{{A}}\in \mathbb{Q}^{mp\times np}$.
\begin{example}\label{exa:0}
Let ${\mathbf{F_p}}$ be the $p\times p$ {\rm DFT} matrix defined by \eqref{eq:DFT}
and ${\mathbf{A}}$ be an arbitrary circulant quaternion matrix ${\mathbf{A}}\in \mathbb{Q}^{p\times p}$ with its unique CD form being $\mathbf{A}=\mathbf{A_1}+\mathbf{A_2}\mathbf{j}$, where $\mathbf{A_1},\mathbf{A_2}\in\mathbb{C}^{p\times p}$.
\begin{itemize}
\item [$\clubsuit$] {\bf First, we show that ${\mathbf{A}}$ cannot be diagonalized by the DFT matrix ${\mathbf{F_p}}$, i.e., ${\mathbf{F_p}}{\mathbf{A}}{\mathbf{F^*_p}}$ is not a diagonal matrix.}

Since $\mathbf{j}c = \bar{c}\mathbf{j}$ for any complex number $c$, it follows that
$\mathbf{j}\mathbf{F_p^*}=\mathbf{F_p}\mathbf{j}.$ Thus, it is deduced that
$$ {\mathbf{F_p}}\mathbf{{A}}{\mathbf{F^*_p}} = {\mathbf{F_p}}(\mathbf{{A}_1}+\mathbf{{A}_2}\mathbf{j}){\mathbf{F^*_p}}
= {\mathbf{F_p}}\mathbf{{A}_1}{\mathbf{F^*_p}}+({\mathbf{F_p}}\mathbf{{A}_2}{\mathbf{F_p}})\mathbf{j}.
$$
As shown in Section \ref{Sect.2}, the first item ``${\mathbf{F_p}}\mathbf{{A}_1}{\mathbf{F^*_p}}$" is a complex diagonal matrix. However, the second item  ``${\mathbf{F_p}}\mathbf{{A}_2}{\mathbf{F_p}}$" is not diagonal in general, as \eqref{eq:F4} shows. Thereby, ${\mathbf{F_p}}\mathbf{{A}}{\mathbf{F^*_p}}$ is not diagonal in general.

\item[$\clubsuit$] {\bf Second, we show that $\mathbf{{A}}$ cannot be diagonalized by two unitary quaternion matrices: ${\mathbf{F_p}}\mathbf{j}$ and $({\mathbf{F_p}}+{\mathbf{F_p}}\mathbf{j})/\sqrt{2}$, respectively. That is to say, neither $({\mathbf{F_p}}\mathbf{j})\mathbf{{A}}({\mathbf{F_p}}\mathbf{j})^*$ nor $\frac{1}{2}({\mathbf{F_p}}+{\mathbf{F_p}}\mathbf{j})\mathbf{{A}}({\mathbf{F_p}}+{\mathbf{F_p}}\mathbf{j})^*$ is a diagonal matrix.}

\begin{itemize}
\item Since $({\mathbf{F_p}}\mathbf{j})({\mathbf{F_p}}\mathbf{j})^*=
    ({\mathbf{F_p}}\mathbf{j})(-{\mathbf{F_p}}\mathbf{j})=\mathbf{I_p}$, we have that ${\mathbf{F_p}}\mathbf{j}$ is unitary.
Further, it is deduced
 $$ ({\mathbf{F_p}}\mathbf{j})\mathbf{{A}}({\mathbf{F_p}}\mathbf{j})^* = ({\mathbf{F_p}}\mathbf{j})(\mathbf{{A}_1}+\mathbf{{A}_2}\mathbf{j})(-{\mathbf{F_p}}\mathbf{j})
= {\mathbf{F_p}}\overline{\mathbf{{A}_1}}{\mathbf{F^*_p}}+({\mathbf{F_p}}\overline{\mathbf{{A}_2}}{\mathbf{F_p}})\mathbf{j},
$$
thereby, ${\mathbf{F_p}\mathbf{j}}\mathbf{{A}}(\mathbf{F_p}\mathbf{j})^*$ is not diagonal in general.

\item Since $ ({\mathbf{F_p}}+{\mathbf{F_p}}\mathbf{j})({\mathbf{F_p}}+{\mathbf{F_p}}\mathbf{j})^*
    =({\mathbf{F_p}}+{\mathbf{F_p}}\mathbf{j})({\mathbf{F^*_p}}-{\mathbf{F_p}}\mathbf{j}) ={\mathbf{F_p}}{\mathbf{F^*_p}}-
{\mathbf{F_p}}{\mathbf{F_p}}\mathbf{j}+{\mathbf{F_p}}{\mathbf{F_p}}\mathbf{j}+{\mathbf{F_p}}{\mathbf{F^*_p}}
=2\mathbf{I_p},$
it follows that $({\mathbf{F_p}}+{\mathbf{F_p}}\mathbf{j})/\sqrt{2}$ is unitary. Further, it is deduced
 $$\begin{array}{rcl}
& & [\frac{1}{\sqrt{2}}({\mathbf{F_p}}+{\mathbf{F_p}}\mathbf{j})]\mathbf{{A}}
[\frac{1}{\sqrt{2}}({\mathbf{F_p}}+{\mathbf{F_p}}\mathbf{j})]^*\\
&=& \frac{1}{2}({\mathbf{F_p}}+{\mathbf{F_p}}\mathbf{j})(\mathbf{{A}_1}
+\mathbf{{A}_2}\mathbf{j})({\mathbf{F^*_p}}-{\mathbf{F_p}}\mathbf{j})\\
&=& \frac{1}{2}[{\mathbf{F_p}}(\mathbf{{A}_1}+\overline{\mathbf{{A}_1}}
+\mathbf{{A}_2}-\overline{\mathbf{{A}_2}})
{\mathbf{F^*_p}}+
{\mathbf{F_p}}(\overline{\mathbf{{A}_1}}-\mathbf{{A}_1}+\mathbf{{A}_2}
+\overline{\mathbf{{A}_2}}){\mathbf{F_p}}\mathbf{j}].
\end{array}$$
thereby, $[\frac{1}{\sqrt{2}}({\mathbf{F_p}}+{\mathbf{F_p}}\mathbf{j})]\mathbf{{A}}
[\frac{1}{\sqrt{2}}({\mathbf{F_p}}+{\mathbf{F_p}}\mathbf{j})]^*$ is not diagonal in general.\end{itemize}
\end{itemize}
In a word, neither unitary complex matrix $\mathbf{F_p}$, nor unitary quaternion matrices $\mathbf{F_p}\mathbf{j}$ and $\frac{1}{\sqrt{2}}(\mathbf{F_p}+\mathbf{F_p}\mathbf{j})$ is a diagonalization matrix of circulant quaternion matrices $\mathbf{A}\in\mathbb{Q}^{p\times p}$.
\end{example}

Example \ref{exa:0} gives negative answers to {\bf Question 1}. Then {\bf Question 2} naturally arises. That is, what kind of unitary quaternion matrices can diagonalize an arbitrary block circulant quaternion matrix? To answer this question, we give the following theorem, which provides sufficient and necessary conditions on unitary quaternion matrix $\mathbf{F}\in \mathbb{Q}^{p\times p}$ being a diagonalization matrix of $p\times p$ circulant quaternion matrices. If $\mathbf{F}\in \mathbb{Q}^{p\times p}$ is a diagonalization matrix of $p\times p$ circulant quaternion matrices, then it can be proved by a similar way as Theorem \ref{theo:block-diagonal} that any given block circulant quaternion matrix of size $mp\times np$ can be block diagonalized by $\mathbf{F}\otimes \mathbf{I_m}\in \mathbb{Q}^{mp\times mp}$ and $\mathbf{F}\otimes \mathbf{I_n}\in \mathbb{Q}^{np\times np}$.
\begin{theorem}\label{theo:1}
Let $\mathbf{{F}}\in \mathbb{Q}^{p\times p}$ be a unitary matrix with the CD form $\mathbf{F}=\mathbf{{F_1}}+\mathbf{{F_2}}\mathbf{j}$, where $\mathbf{\mathbf{F_1}}$, $\mathbf{\mathbf{F_2}}\in \mathbb{C}^{p\times p}$.
Then $\mathbf{F}$ is a diagonalization matrix of $p\times p$ circulant quaternion matrices if and only if all these matrices
\begin{eqnarray}\label{condition-1}
\mathbf{\mathbf{F_1}}\mathbf{\mathbf{\hat{A}}}{\mathbf{{F^*_1}}},\  \mathbf{\mathbf{F_1}}\mathbf{\mathbf{\hat{A}}}\mathbf{{F^\top_1}},\ \mathbf{\mathbf{F_1}}\mathbf{\mathbf{\hat{A}}}{\mathbf{{F}^\top_2}},\ \mathbf{F_1}\mathbf{\hat{A}}\mathbf{F_2^*},\nonumber \\
\mathbf{F_2}\mathbf{\hat{A}}{{\mathbf{F_1^*}}},\ \mathbf{F_2}\mathbf{\hat{A}}{\mathbf{F_1^\top}},\ \mathbf{F_2}\mathbf{\hat{A}}{\mathbf{F_2^\top}},\ \mathbf{F_2}\mathbf{\hat{A}}\mathbf{F_2^*},
\end{eqnarray}
are diagonal for any given real circulant matrices $\mathbf{\hat{A}}=\operatorname{circ}(\hat{a}_1,\hat{a}_2,\ldots, \hat{a}_p)\in\mathbb{R}^{p\times p}$.
\end{theorem}

\noindent{\bf Proof.} For any circulant quaternion matrix $\mathbf{{A}}=\operatorname{circ}(a_1,a_2,\ldots, a_p)\in \mathbb{Q}^{p\times p}$, rewrite it by the CD form  $\mathbf{A}=\mathbf{{A_1}}+\mathbf{{A_2}}\mathbf{j}$ with $\mathbf{{A_1}},\mathbf{{A_2}}\in \mathbb{C}^{p\times p}$. Since $\mathbf{{F}^*}=(\mathbf{F_1^*})-(\mathbf{F_2^\top})\mathbf{j}$ from \eqref{conj}, it can be deduced
\begin{equation}\label{0}
\begin{array}{rcl}
\mathbf{{F}}\mathbf{{A}}\mathbf{{F}^*}
&=&(\mathbf{F_1}+\mathbf{F_2}\mathbf{j})(\mathbf{{A}_1}+\mathbf{{A}_2}\mathbf{j})
(\mathbf{F_1^*}-\mathbf{F_2^\top}\mathbf{j})\\
&=& (\mathbf{F_1}\mathbf{{A}_1}\mathbf{F_1^*}-\mathbf{F_2}\overline{\mathbf{{A}_2}}\mathbf{F_1^*}+
\mathbf{F_2}\overline{\mathbf{{A}_1}}\mathbf{F_2^*}+\mathbf{F_1}\mathbf{{A}_2}\mathbf{F_2^*})
\\
&&-
(\mathbf{F_1}\mathbf{{A}_1}\mathbf{F_2^\top}+\mathbf{F_2}\overline{\mathbf{{A}_2}}\mathbf{F_2^\top}+
\mathbf{F_2}\overline{\mathbf{{A}_1}}\mathbf{F_1^\top}-\mathbf{F_1}\mathbf{{A}_2}\mathbf{F_1^\top})\mathbf{j}
\end{array}
\end{equation}
Hence, the sufficiency is obvious. Below, we prove the necessity.

Assume that $\mathbf{\hat{A}}=\operatorname{circ}(\hat{a}_1,\hat{a}_2,\ldots, \hat{a}_p)$ is an any given circulant real matrix in $\mathbb{R}^{p\times p}$ and $\mathbf{{F}}\mathbf{{A}}\mathbf{{F}^*}$ is a diagonal matrix for any circulant quaternion matrix $\mathbf{A}$.
Then we take $\mathbf{{A}}$ as the following special cases and substitute it into \eqref{0} to deduce the following results:
\begin{itemize}
\item Take $\mathbf{{A}_1}=\mathbf{{A}_2}=\mathbf{\hat{A}}$. Then from \eqref{0} and $\mathbf{{F}}\mathbf{{A}}\mathbf{{F}^*}$ being diagonal, it follows that
\begin{equation}\label{1}
\left\{\begin{array}{ll}
\mathbf{F_1}\mathbf{\hat{A}}\mathbf{F_1^*}-\mathbf{F_2}\mathbf{\hat{A}}\mathbf{F_1^*}-
\mathbf{F_2}\mathbf{\hat{A}}(-\mathbf{F_2^*})-\mathbf{F_1}\mathbf{\hat{A}}(-\mathbf{F_2^*})
&{\mbox{\rm is diagonal}}; \\
\mathbf{F_1}\mathbf{\hat{A}}(-\mathbf{F_2^\top})-\mathbf{F_2}\mathbf{\hat{A}}(-\mathbf{F_2^\top})+
\mathbf{F_2}\mathbf{\hat{A}}\mathbf{F_1^\top}-\mathbf{F_1}\mathbf{\hat{A}}\mathbf{F_1^\top}& {\mbox{\rm is diagonal}}.
\end{array}\right.
\end{equation}
\item Take $\mathbf{{A}_1}=\mathbf{\hat{A}}\mathbf{i}$ and $\mathbf{{A}_2}=-\mathbf{\hat{A}}$. Then from \eqref{0} and $\mathbf{{F}}\mathbf{{A}}\mathbf{{F}^*}$  being diagonal, it follows that
$$\left\{
\begin{array}{ll}
[\mathbf{F_1}\mathbf{\hat{A}}\mathbf{F_1^*}+
\mathbf{F_2}\mathbf{\hat{A}}(-\mathbf{F_2^*})]\mathbf{i}
+[\mathbf{F_2}\mathbf{\hat{A}}\mathbf{F_1^*}+
\mathbf{F_1}\mathbf{\hat{A}}(-\mathbf{F_2^*})]&{\mbox{\rm is diagonal}};\\

[\mathbf{F_1}\mathbf{\hat{A}}(-\mathbf{F_2^\top})-\mathbf{F_2}\mathbf{\hat{A}}\mathbf{F_1^\top}]\mathbf{i}+
[\mathbf{F_2}\mathbf{\hat{A}}(-\mathbf{F_2^\top})
+\mathbf{F_1}\mathbf{\hat{A}}\mathbf{F_1^\top}] &{\mbox{\rm is diagonal}},
\end{array}\right.
$$
which, together with \eqref{1}, imply that
\begin{equation}\label{2}
\left\{\begin{array}{ll}
[\mathbf{F_1}\mathbf{\hat{A}}\mathbf{F_1^*}-
\mathbf{F_2}\mathbf{\hat{A}}(-\mathbf{F_2^*})]+[\mathbf{F_1}\mathbf{\hat{A}}\mathbf{F_1^*}+
\mathbf{F_2}\mathbf{\hat{A}}(-\mathbf{F_2^*})]\mathbf{i}&{\mbox{\rm is diagonal}};\\

[\mathbf{F_1}\mathbf{\hat{A}}(-\mathbf{F_2^\top})+
\mathbf{F_2}\mathbf{\hat{A}}\mathbf{F_1^\top}]+[\mathbf{F_1}\mathbf{\hat{A}}(-\mathbf{F_2^\top})-
\mathbf{F_2}\mathbf{\hat{A}}\mathbf{F_1^\top}]\mathbf{i}&{\mbox{\rm is diagonal}}.
\end{array}\right.
\end{equation}

\item Take $\mathbf{{A}_2}=\mathbf{0}$ and $\mathbf{{A}_1}=\mathbf{\hat{A}}$. Then from \eqref{0} and $\mathbf{{F}{A}{F}^*}$ being diagonal, it follows that
\begin{equation}\label{3}
\left\{\begin{array}{ll}
\mathbf{F_1}\mathbf{\hat{A}}\mathbf{F_1^*}-
\mathbf{F_2}{\mathbf{\hat{A}}}(-\mathbf{F_2^*})&{\mbox{\rm is diagonal}}; \\
\mathbf{F_1}\mathbf{\hat{A}}(-\mathbf{F_2^\top})+
\mathbf{F_2}\mathbf{\hat{A}}\mathbf{F_1^\top}&{\mbox{\rm is diagonal}},
\end{array}\right.
\end{equation}
which together with \eqref{2} imply that
\begin{equation}\label{4}
\left\{\begin{array}{ll}
\mathbf{F_1}\mathbf{\hat{A}}\mathbf{F_1^*}+
\mathbf{F_2}\mathbf{\hat{A}}(-\mathbf{F_2^*}) & {\mbox{\rm is diagonal}}; \\
\mathbf{F_1}\mathbf{\hat{A}}(-\mathbf{F_2^\top})-
\mathbf{F_2}\mathbf{\hat{A}}\mathbf{F_1^\top} & {\mbox{\rm is diagonal}}.
\end{array}\right.
\end{equation}
By equations \eqref{3} and \eqref{4}, it further follows that
\begin{equation}\label{5}
\mathbf{F_1}\mathbf{\hat{A}}\mathbf{F_1^*},\
\mathbf{F_2}\mathbf{\hat{A}}\mathbf{F_2^*},\
\mathbf{F_1}\mathbf{\hat{A}}\mathbf{F_2^\top}\  {\mbox{\rm and}}\
\mathbf{F_2}\mathbf{\hat{A}}\mathbf{F_1^\top},\quad {\mbox{\rm are diagonal}}.
\end{equation}

\item Take $\mathbf{{A}_2}=\mathbf{\hat{A}}\mathbf{i}$ and $\mathbf{{A}_1}=-\mathbf{\hat{A}}$. Then from \eqref{0} and $\mathbf{{F}{A}{F}^*}$     being diagonal, it follows that
\begin{equation*}
\left\{\begin{array}{ll}
[\mathbf{F_2}\mathbf{\hat{A}}\mathbf{F_1^*}-
\mathbf{F_1}\mathbf{\hat{A}}(-\mathbf{F_2^*})]\mathbf{i}+[-\mathbf{F_1}\mathbf{\hat{A}}\mathbf{F_1^*}+
\mathbf{F_2}\mathbf{\hat{A}}(-\mathbf{F_2^*})]&{\mbox{\rm is diagonal}};\\

[\mathbf{F_2}\mathbf{\hat{A}}(-\mathbf{F_2^\top})-
\mathbf{F_1}\mathbf{\hat{A}}\mathbf{F_1^\top}]\mathbf{i}+[-\mathbf{F_1}\mathbf{\hat{A}}(-\mathbf{F_2^\top})-
\mathbf{F_2}\mathbf{\hat{A}}(F_1)^\top]&{\mbox{\rm is diagonal}},
\end{array}\right.
\end{equation*}

which, together with \eqref{1}, imply that
\begin{equation}\label{6}
\left\{\begin{array}{ll}
[\mathbf{F_2}\mathbf{\hat{A}}\mathbf{F_1^*}-
\mathbf{F_1}\mathbf{\hat{A}}(-\mathbf{F_2^*})]\mathbf{i}-[\mathbf{F_2}\mathbf{\hat{A}}\mathbf{F_1^*}+
\mathbf{F_1}\mathbf{\hat{A}}(-\mathbf{F_2^*})] & {\mbox{\rm is diagonal}};\\

[\mathbf{F_2}\mathbf{\hat{A}}(-\mathbf{F_2^\top})-
\mathbf{F_1}\mathbf{\hat{A}}\mathbf{F_1^\top}]\mathbf{i}-[\mathbf{F_2}\mathbf{\hat{A}}\mathbf{(-F_2^\top)}+
\mathbf{F_1}\mathbf{\hat{A}}\mathbf{F_1^\top}]&{\mbox{\rm is diagonal}}.
\end{array}\right.
\end{equation}

\item Take $\mathbf{{A}_1}=\mathbf{0}$ and $\mathbf{{A}_2}=\mathbf{\hat{A}}$. Then from \eqref{0} and $\mathbf{{F}{A}{F}^*}$     being diagonal, it follows that
\begin{equation}\label{7}
\left\{\begin{array}{ll}
\mathbf{F_2}\mathbf{\hat{A}}\mathbf{F_1^*}+
\mathbf{F_1}{\mathbf{\hat{A}}}(-\mathbf{F_2^*})&{\mbox{\rm is diagonal}}; \\

\mathbf{F_2}\mathbf{\hat{A}}(-\mathbf{F_2^\top})+
\mathbf{F_1}\mathbf{\hat{A}}\mathbf{F_1^\top}&{\mbox{\rm is diagonal}},
\end{array}\right.
\end{equation}
which, together with \eqref{6}, imply that
\begin{equation}\label{8}
\left\{\begin{array}{ll}
\mathbf{F_2}\mathbf{\hat{A}}\mathbf{F_1^*}-
\mathbf{F_1}\mathbf{\hat{A}}(-\mathbf{F_2^*})&{\mbox{\rm is diagonal}}; \\

\mathbf{F_2}\mathbf{\hat{A}}(-\mathbf{F_2^\top})-
\mathbf{F_1}\mathbf{\hat{A}}\mathbf{F_1^\top}& {\mbox{\rm is diagonal}},
\end{array}\right.
\end{equation}
By equations \eqref{7} and \eqref{8}, it further follows that
\begin{equation}\label{9}
\mathbf{F_2}\mathbf{\hat{A}}\mathbf{F_1^*},\
\mathbf{F_1}\mathbf{\hat{A}}\mathbf{F_2^*},\
\mathbf{F_2}\mathbf{\hat{A}}\mathbf{F_2^\top}\ {\mbox{\rm and}}\
\mathbf{F_1}\mathbf{\hat{A}}\mathbf{F_1^\top},\quad  {\mbox{\rm are diagonal}}.
\end{equation}
\end{itemize}

By combining equations \eqref{5} and \eqref{9}, the necessity is obtained. The proof is completed. \ep

\begin{remark}
${\rm (i)}$  Theorem \ref{theo:1} provides an alternative way to find unitary quaternion diagonalization matrix $\mathbf{F}=\mathbf{F_1}+\mathbf{F_2}\mathbf{j}$ for circulant quaternion matrix, that is to find complex matrices $\mathbf{F_1}$ and $\mathbf{F_2}$ such that $\mathbf{F}$ is unitary and all these matrices in \eqref{condition-1} are diagonal for any given real circulant matrix $\mathbf{\mathbf{\hat{A}}}$. However, to find such complex matrices $\mathbf{F_1}$ and $\mathbf{F_2}$ is too hard. In addition, Example \ref{exa:0} has demonstrated that both $\mathbf{F_1}$ and $\mathbf{F_2}$ in Theorem \ref{theo:1} cannot be the DFT matrix, which implies that even if such complex matrices $\mathbf{F_1}$ and $\mathbf{F_2}$ that satisfy conditions in Theorem \ref{theo:1} are found, we cannot make use of the fast computation way by performing FFT on complex vectors to obtain the corresponding result diagonal matrix, which is not what we want. Hence, we are no longer committed to finding unitary diagonalization matrices of circulant quaternion matrices in the quaternion domain.

${\rm (ii)}$ To establish a way to fast obtain the result diagonal quaternion matrix of a circulant quaternion matrix by making use of the calculation convenience of FFTs, we should turn to a larger domain than the quaternion domain to achieve the diagonalization of circulant quaternion matrices. From the discussion of Example \ref{exa:0}, it can be seen that the appearance of this part ``$\mathbf{F_p}\mathbf{\hat{A}}\mathbf{F_p}$" where $\mathbf{\hat{A}}$ is a complex matrix, is the main reason that leads to that $(\mathbf{F_p}\mathbf{p})\mathbf{A}(\mathbf{F_p}\mathbf{p})^*$ with $\mathbf{p}=1,\mathbf{j}$ or $(1+\mathbf{j})/\sqrt{2}$ cannot be diagonal for any given circulant quaternion matrix $\mathbf{A}$. The appearance of this part ``$\mathbf{F_p}\mathbf{\hat{A}}\mathbf{F_p}$" is due to the special commutative cases of the quaternion algebra that $1c=c1$ but $\mathbf{j}c=\bar{c}\mathbf{j}$ for any complex number $c$. Thus, noting that $\mathbf{l}q=\bar{q}\mathbf{l}$ for any quaternion $q$, which helps to avoid inconsistent commutative cases of $1c=c1$ and $\mathbf{j}c=\bar{c}\mathbf{j}$, we study the diagonalization of circulant quaternion matrices and further the block diagonalization of block circulant quaternion matrices in the octonion domain, in the next section.
\end{remark}

\section{Diagonalization matrices in the octonion domain}
\label{Sect.3.2}
In this section, we aim to find a unitary octonion matrix $\mathbf{\mathbf{O}}\in \mathbb{O}^{p\times p}$ such that
$$(\mathbf{O}\otimes \mathbf{I_m})\mathbf{Q}(\mathbf{O}\otimes \mathbf{I_n})^*={\rm Diag}(\mathbf{\hat{Q}_1},\mathbf{\hat{Q}_2},\ldots,\mathbf{\hat{Q}_p})\in \mathbb{Q}^{mp\times np}$$
for any block circulant quaternion matrix $\mathbf{Q}=\operatorname{bcirc}(\mathbf{Q_1},\mathbf{Q_2},\ldots,\mathbf{Q_p})\in \mathbb{Q}^{p\times p}$. The octonion algebra has some wonderful properties what we want, for instance, $\mathbf{l}q=\bar{q}\mathbf{l}$ for any quaternion number $q$ and $(\mathbf{sl})(\mathbf{tl})\in \mathbb{Q}$ for any $s,t\in \mathbb{Q}$, but at the same time the octonion algebra bring some undesirable facts, that is the octonion algebra is neither commutative nor associative. Therefore, in Subsection \ref{Sect.3-proposition} we give several special cases of the octonion multiplication, based on which we achieve the diagonalization of circulant quaternion matrices and the block diagonalization of block circulant quaternion matrices in Subsection \ref{subsect.diag-oct} and  Subsection \ref{subsect.blockdiag-oct}, respectively.
\subsection{Several special cases of the octonion multiplication}
\label{Sect.3-proposition}
\setcounter{equation}{0} \setcounter{assumption}{0}
\setcounter{theorem}{0} \setcounter{proposition}{0}
\setcounter{corollary}{0} \setcounter{lemma}{0}
\setcounter{definition}{0} \setcounter{remark}{0}
\setcounter{algorithm}{0}
From Table \ref{multiplication-octonion}, it can be easily seen that the octonion multiplication is non-commutative and non-associative, for example:
$$\mathbf{i}\mathbf{j}\neq\mathbf{j}\mathbf{i},\quad{\rm and}\quad \mathbf{i}(\mathbf{jl})\neq (\mathbf{ij})\mathbf{l}.$$
This undesirable nature makes the octonion multiplication very complex and difficult to use. For the follow-up study, we give the following two propositions on special commutative cases and associative cases of the octonion multiplication, respectively.

\begin{proposition}\label{pro:octonion}
Let $p=a+b\mathbf{i}$, $q=c+d\mathbf{i}$, $r=e+f\mathbf{i}$ and $s=g+h\mathbf{i}$ with $a,b,c,d,e,f,g,h\in\mathbb{R}$. Then
\begin{itemize}
\item [{\mbox{\rm (i)}}]
    $(p+q\mathbf{j})\cdot[(r+s\mathbf{j})\mathbf{l}]=[(r+s\mathbf{j})\mathbf{l}] \cdot\overline{(p+q\mathbf{j})},$
\item [{\mbox{\rm (ii)}}] $(p\mathbf{l})\cdot q\cdot (r\mathbf{l})=-p\bar{q}\bar{r}=[p(\mathbf{jl})]\cdot q\cdot [r(\mathbf{jl})]$,
\item [{\mbox{\rm (iii)}}]
$(p\mathbf{l})\cdot (q\mathbf{j})\cdot(r\mathbf{l})=pq\bar{r}\mathbf{j}=[p(\mathbf{jl})]\cdot (q\mathbf{j})\cdot[r(\mathbf{jl})]$,

\item [{\mbox{\rm (iv)}}]
$(p\mathbf{l})\cdot q\cdot[r(\mathbf{jl})]=\bar{p}q\bar{r}\mathbf{j}=-[p(\mathbf{jl})]\cdot q\cdot(r\mathbf{l})$,
\item [{\mbox{\rm (v)}}]
$(p\mathbf{l})\cdot (q\mathbf{j})\cdot[r(\mathbf{jl})]=\bar{p}\bar{q}\bar{r}=-[p(\mathbf{jl})]\cdot (q\mathbf{j})\cdot(r\mathbf{l})$.
\end{itemize}
\end{proposition}

{\bf Proof.} {\mbox{\rm (i)}} From the octonion multiplication rules in Table \ref{multiplication-octonion}, it can be seen that
if $u\in\{\mathbf{l},\mathbf{il},\mathbf{jl},\mathbf{kl}\}$, then $(p+q\mathbf{j})u=u\overline{(p+q\mathbf{j})}=u(\bar{p}-q\mathbf{j}).$
Thus, it can be deduced
$$(p+q\mathbf{j})\cdot[(r+s\mathbf{j})\mathbf{l}] =[(r+s\mathbf{j})\mathbf{l}]\cdot\overline{(p+q\mathbf{j})}.$$

{\mbox{\rm (ii)}} From the octonion multiplication rules in Table \ref{multiplication-octonion}, it follows that
$$\begin{array}{rcl}
(p\mathbf{l})\cdot q\cdot(r\mathbf{l})
&=&(a\mathbf{l}+b\mathbf{il})(c+d\mathbf{i})(e\mathbf{l}+f\mathbf{il})\\
&=&-ace-bce\mathbf{i}+ade\mathbf{i}-bde
+acf\mathbf{i}-bcf+adf+bdf\mathbf{i}\\
&=&-(a+b\mathbf{i})(c-d\mathbf{i})(e-f\mathbf{i})=-p\bar{q}\bar{r},
\end{array}$$
$$\begin{array}{rcl}
[p(\mathbf{jl})]\cdot q\cdot[r(\mathbf{jl})]&=&(a\mathbf{jl}+b\mathbf{kl})(c+d\mathbf{i})(e\mathbf{jl}+f\mathbf{kl})\\
&=&-ace-bce\mathbf{i}+ade\mathbf{i}-bde
+acf\mathbf{i}-bcf+adf+bdf\mathbf{i}\\
&=&-(a+b\mathbf{i})(c-d\mathbf{i})(e-f\mathbf{i})=-p\bar{q}\bar{r}.
\end{array}$$
Thus, it can be deduced $(p\mathbf{l})\cdot q\cdot(r\mathbf{l})=-p\bar{q}\bar{r}=[p(\mathbf{jl})]\cdot q\cdot [r(\mathbf{jl})].$

{\mbox{\rm (iii)}} Similar to the proof of {\mbox{\rm (ii)}}, from Table \ref{multiplication-octonion}, we can obtain that
$$\begin{array}{rcl}
(p\mathbf{l})\cdot(q\mathbf{j})\cdot(r\mathbf{l})
&=&(a\mathbf{l}+b\mathbf{il})(c\mathbf{j}+d\mathbf{k})(e\mathbf{l}+f\mathbf{il})\\
&=&ace\mathbf{j}+bce\mathbf{k}+
ade\mathbf{k}-bde\mathbf{j}-acf\mathbf{k}+bcf\mathbf{j}+
adf\mathbf{j}+bdf\mathbf{k}\\
&=&(a+b\mathbf{i})(c+d\mathbf{i})(e-f\mathbf{i})\mathbf{j}=pq\bar{r}\mathbf{j},
\end{array}$$
$$\begin{array}{rcl}
[p(\mathbf{jl})]\cdot (q\mathbf{j})\cdot [r(\mathbf{jl})]&=& (a\mathbf{jl}+b\mathbf{kl})(c\mathbf{j}+d\mathbf{k})(e\mathbf{jl}+f\mathbf{kl})\\
&=&ace\mathbf{j}+bce\mathbf{k}+
ade\mathbf{k}-bde\mathbf{j}-acf\mathbf{k}+bcf\mathbf{j}+
adf\mathbf{j}+bdf\mathbf{k}\\
&=&(a+b\mathbf{i})(c+d\mathbf{i})(e-f\mathbf{i})\mathbf{j}=pq\bar{r}\mathbf{j},
\end{array}$$
thus, it follows that
$(p\mathbf{l})\cdot (q\mathbf{j})\cdot(r\mathbf{l})=pq\bar{r}\mathbf{j}=[p(\mathbf{jl})]\cdot (q\mathbf{j})\cdot[r(\mathbf{jl})].$

{\mbox{\rm (iv)}} Similarly, from Table \ref{multiplication-octonion}, we can obtain that
$$\begin{array}{rcl}
(p\mathbf{l})\cdot q\cdot[r(\mathbf{jl})]
&=&(a\mathbf{l}+b\mathbf{il})(c+d\mathbf{i})(e\mathbf{jl}+f\mathbf{kl})\\
&=&(ace\mathbf{j}-bce\mathbf{k}+ade\mathbf{k}+bde\mathbf{j}
-acf\mathbf{k}-bcf\mathbf{j}+adf\mathbf{j}-bdf\mathbf{k}\\
&=&(a-b\mathbf{i})(c+d\mathbf{i})(e-f\mathbf{i})\mathbf{j}=\bar{p}{q}\bar{r}\mathbf{j},
\end{array}$$
$$\begin{array}{rcl}
[p(\mathbf{jl})]\cdot q\cdot(r\mathbf{l})
&=&(a\mathbf{jl}+b\mathbf{kl})(c+d\mathbf{i})(e\mathbf{l}+f\mathbf{il})\\
&=&(ace\mathbf{j}-bce\mathbf{k}+ade\mathbf{k}+bde\mathbf{j}
-acf\mathbf{k}-bcf\mathbf{j}+adf\mathbf{j}-bdf\mathbf{k}\\
&=&(a-b\mathbf{i})(c+d\mathbf{i})(e-f\mathbf{i})\mathbf{j}=\bar{p}{q}\bar{r}\mathbf{j},
\end{array}$$
thus, it follows that
$(p\mathbf{l})\cdot q\cdot[r(\mathbf{jl})]=\bar{p}q\bar{r}\mathbf{j}=-[p(\mathbf{jl})]\cdot q\cdot(r\mathbf{l}).$

{\mbox{\rm (v)}} Similarly, from Table \ref{multiplication-octonion}, we can obtain that
$$\begin{array}{rcl}
(p\mathbf{l})\cdot(q\mathbf{j})\cdot[r(\mathbf{jl})]
&=&(a\mathbf{l}+b\mathbf{il})(c\mathbf{j}+d\mathbf{k})(e\mathbf{jl}+f\mathbf{kl})\\
&=&ace-bce\mathbf{i}-ade\mathbf{i}-bde
-acf\mathbf{i}-bcf-adf+bdf\mathbf{i}\\
&=&(a-b\mathbf{i})(c-d\mathbf{i})(e-f\mathbf{i})=\bar{p}\bar{q}\bar{r},
\end{array}$$
$$\begin{array}{rcl}
[p(\mathbf{jl})]\cdot(q\mathbf{j})\cdot(r\mathbf{l})
&=&(a\mathbf{jl}+b\mathbf{kl})(c\mathbf{j}+d\mathbf{k})(e\mathbf{l}+f\mathbf{il})\\
&=&ace-bce\mathbf{i}-ade\mathbf{i}-bde
-acf\mathbf{i}-bcf-adf+bdf\mathbf{i}\\
&=&(a-b\mathbf{i})(c-d\mathbf{i})(e-f\mathbf{i})=\bar{p}\bar{q}\bar{r},
\end{array}$$
thus, it follows that
$(p\mathbf{l})\cdot (q\mathbf{j})\cdot[r(\mathbf{jl})]=\bar{p}\bar{q}\bar{r}=-[p(\mathbf{jl})]\cdot (q\mathbf{j})\cdot(r\mathbf{l}).$\ep

\begin{remark}
\mbox{\rm (i)} From Proposition \ref{pro:octonion}, it can be seen that for any two quaternions $p$ and $q$, $pq\neq qp$ in general, but they satisfy a special commutative rule, i.e., $p(q\mathbf{l})=(q\mathbf{l})\bar{p}.$
In addition, it should be noticed that
$p(q\mathbf{l})\neq\overline{(q\mathbf{l})}{p}$
 in general.

\mbox{\rm (ii)} Let $\mathbf{p}=\mathbf{l},\mathbf{jl}$ or $(\mathbf{l}+\mathbf{jl})/\sqrt{2}$. By {\rm (ii)}, {\rm (iii)}, {\rm (iv)} and {\rm (v)} of Proposition \ref{pro:octonion}, it follows that for any $p, q, r, s\in \mathbb{C}$,
\begin{equation}\label{eq:relation-remark}
  (p\mathbf{p})\cdot (q+r\mathbf{j})\cdot(s\mathbf{p})=-p\bar{q}\bar{s}+pr\bar{s}\mathbf{j}.
\end{equation}
This fact will be frequently used in the follow-up discussion.
\end{remark}

The next proposition is about special associative cases of the octonion multiplication.
\begin{proposition}\label{pro:associative}
Let $p=a+b\mathbf{i}$, $q=c+d\mathbf{i}$, $r=e+f\mathbf{i}$, $s=g+h\mathbf{i}$ and $t=m+n\mathbf{i}$ with $a,b,c,d,e,f,g,h, m,n\in\mathbb{R}$.
Then
\begin{itemize}
\item [{\mbox{\rm (i)}}] $p[(q\mathbf{l})r(s\mathbf{l})]=p(q\mathbf{l})r(s\mathbf{l})$, $p[(q\mathbf{l})(r\mathbf{j})(s\mathbf{l})]=p(q\mathbf{l})(r\mathbf{j})(s\mathbf{l})$.
\item [{\mbox{\rm (ii)}}] $(p\mathbf{j})[(q\mathbf{l})r(s\mathbf{l})]=(p\mathbf{j})(\bar{q}\mathbf{l})r(\bar{s}\mathbf{l})$,
$(p\mathbf{j})[(q\mathbf{l})(r\mathbf{j})(s\mathbf{l})]=(p\mathbf{j})(\bar{q}\mathbf{l})(r\mathbf{j})(\bar{s}\mathbf{l})$.
\item [{\mbox{\rm (iii)}}] $(p\mathbf{l})[(q\mathbf{l})r(s\mathbf{l})](t\mathbf{l})
    =[(p\mathbf{l})(q\mathbf{l})]r[(s\mathbf{l})(t\mathbf{l})]$,
    $(p\mathbf{l})[qrs](t\mathbf{l})
    =[(p\mathbf{l})q]r[s(t\mathbf{l})]$.

\item [\mbox{\rm (iv)}] $(p\mathbf{l})[(q\mathbf{l})(r\mathbf{j})(s\mathbf{l})](t\mathbf{l})
    =[({p}\mathbf{l})(\bar{q}\mathbf{l})](r\mathbf{j})[(s\mathbf{l})(\bar{t}\mathbf{l})]$, $(p\mathbf{l})[c(r\mathbf{j})d](s\mathbf{l})
    =[({p}\mathbf{l})c](r\mathbf{j})[d({s}\mathbf{l})]$. 

\item [\mbox{\rm (v)}] $(p\mathbf{l})(q\mathbf{l})(r\mathbf{l})=(p\mathbf{l})[(q\mathbf{l})(r\mathbf{l})]$,
$(p\mathbf{jl})(q\mathbf{l})(r\mathbf{l})=(p\mathbf{jl})[(\bar{q}\mathbf{l})(\bar{r}\mathbf{l})]$.

\item [\mbox{\rm (vi)}] $(t\mathbf{l})(p+q\mathbf{j})(r+s\mathbf{j})=(t\mathbf{l})[pr-\bar{q}s+(\bar{p}s+qr)\mathbf{j}]$.
\end{itemize}
\end{proposition}

{\bf Proof.}
{\mbox{\rm (i)}} From Proposition \ref{pro:octonion}, it follows that
\begin{eqnarray*}
&&p[(q\mathbf{l})r(s\mathbf{l})]=p(-q\bar{r}\bar{s}),\ %
 p(q\mathbf{l})r(s\mathbf{l})=(q\mathbf{l})\bar{p}r(s\mathbf{l})=-qp\bar{r}\bar{s},\\
&&p[(q\mathbf{l})(r\mathbf{j})(s\mathbf{l})]=p(qr\bar{s}\mathbf{j}), \ %
 p(q\mathbf{l})(r\mathbf{j})(s\mathbf{l})=(q\mathbf{l})\bar{p}(r\mathbf{j})(s\mathbf{l})
=qpr\bar{s}\mathbf{j}.
\end{eqnarray*}
Noting that $p,q,r,s\in \mathbb{C}$, thus we have that
$$p[(q\mathbf{l})r(s\mathbf{l})]=p(q\mathbf{l})r(s\mathbf{l})\quad\mbox{\rm and}\quad p[(q\mathbf{l})(r\mathbf{j})(s\mathbf{l})]=p(q\mathbf{l})(r\mathbf{j})(s\mathbf{l}).$$

{\mbox{\rm (ii)}} From Proposition \ref{pro:octonion} (i) and (ii), it follows that $$(p\mathbf{j})[(q\mathbf{l})r(s\mathbf{l})]=-p\bar{q}rs\mathbf{j}.$$
By the octonion multiplication rules shown in Table \ref{multiplication-octonion}, we can obtain that 
$$\begin{array}{rcl}
& &(p\mathbf{j})(q\mathbf{l})(r\mathbf{j})(s\mathbf{l})\\
&=&(a\mathbf{j}+b\mathbf{k})(c\mathbf{l}+d\mathbf{il})(e\mathbf{j}+f\mathbf{k})(g\mathbf{l}+h\mathbf{il})\\
&=&-aceg\mathbf{j}-bceg\mathbf{k}-adeg\mathbf{k}+bdeg\mathbf{j}
-acfg\mathbf{k}+bcfg\mathbf{j}+adfg\mathbf{j}+bdfg\mathbf{k}\\
& &+aceh\mathbf{k}-bceh\mathbf{j}-adeh\mathbf{j}-bdeh\mathbf{k}
-acfh\mathbf{j}-bcfh\mathbf{k}-adfh\mathbf{k}+bdfh\mathbf{j}\\
&=&-pqr\bar{s}\mathbf{j}\\
\end{array}$$
From Proposition \ref{pro:octonion} (i) and (iii), it follows that
$$(p\mathbf{j})[(q\mathbf{l})(r\mathbf{j})(s\mathbf{l})]=(p\mathbf{j})(qr\bar{s}\mathbf{j})=-p\bar{q}\bar{r}s.$$
By the octonion multiplication rules shown in Table \ref{multiplication-octonion}, we can obtain that
$$\begin{array}{rcl}
& &(p\mathbf{j})(q\mathbf{l})(r\mathbf{j})(s\mathbf{l})\\
&=&(a\mathbf{j}+b\mathbf{k})(c\mathbf{l}+d\mathbf{il})(e\mathbf{j}+f\mathbf{k})(g\mathbf{l}+h\mathbf{il})\\
&=&-aceg-bceg\mathbf{i}-adeg\mathbf{i}+bdeg
+acfg\mathbf{i}-bcfg-adfg-bdfg\mathbf{i}\\
& &+aceh\mathbf{i}-bceh-adeh-bdeh\mathbf{i}
+acfh+bcfh\mathbf{i}+adfh\mathbf{i}-bdfh\\
&=&-pq\bar{r}\bar{s}.
\end{array}$$
Hence, it follows that $$(p\mathbf{j})[(q\mathbf{l})r(s\mathbf{l})]=(p\mathbf{j})(\bar{q}\mathbf{l})r(\bar{s}\mathbf{l})\quad\mbox{\rm and}\quad (p\mathbf{j})[(q\mathbf{l})(r\mathbf{j})(s\mathbf{l})]=(p\mathbf{j})(\bar{q}\mathbf{l})(r\mathbf{j}(\bar{s}\mathbf{l}).$$

In addition, by the octonion multiplication rules shown in Table \ref{multiplication-octonion}, it is similar to obtain that {\mbox{\rm (iii)}}, {\mbox{\rm (iv)}}, {\mbox{\rm (v)}} and {\mbox{\rm (vi)}} hold.
\ep

\begin{remark}
It should be noticed from Proposition \ref{pro:associative} that for any $p,q,r, s, t\in \mathbb{C}$,
$$(p\mathbf{l})(q\mathbf{l})(r\mathbf{j})\neq({p}\mathbf{l})[({q}\mathbf{l})(r\mathbf{j})], \;\; (p\mathbf{j})(q\mathbf{l})(r\mathbf{l})\neq(p\mathbf{j})[({q}\mathbf{l})({r}\mathbf{l})],$$
$$ (p\mathbf{jl})(q\mathbf{l})(r\mathbf{l})\neq(p\mathbf{jl})[({q}\mathbf{l})({r}\mathbf{l})],\ \ (t\mathbf{l})(p+q\mathbf{j})(r+s\mathbf{j})\neq (t\mathbf{l})[(p+q\mathbf{j})(r+s\mathbf{j})].$$
\end{remark}


\subsection{Diagonalization of circulant quaternion matrices}
\label{subsect.diag-oct}
In this section, we achieve the diagonalization of circulant quaternion matrices in the octonion domain by making use of Proposition \ref{pro:octonion} established in Section \ref{Sect.3-proposition}.

First, we show that octonion matrices ${\mathbf{F_p}}\mathbf{l}$, ${\mathbf{F_p}}(\mathbf{jl})$ and ${\mathbf{F_p}}(\mathbf{l}+\mathbf{jl})/\sqrt{2}$ are unitary.

\begin{theorem}\label{ther:unitary}
Let ${\mathbf{F_p}}$ be a $p\times p$ DFT matrix defined by \eqref{eq:DFT}. Then octonion matrices ${\mathbf{F_p}}\mathbf{l}$, ${\mathbf{F_p}}(\mathbf{jl})$ and ${\mathbf{F_p}}(\mathbf{l}+\mathbf{jl})/\sqrt{2}$ are unitary.
\end{theorem}

{\bf Proof.} Let $\mathbf{p}=\mathbf{l}, \mathbf{jl}$ or $ (\mathbf{l}+\mathbf{jl})/\sqrt{2}$, and $\mathbf{O}={\mathbf{F_p}}\mathbf{p}$. We show that $\mathbf{OO^*}=\mathbf{I_p}$.

Obviously, $({\mathbf{F_p}}\mathbf{l})^*=-{\mathbf{F_p}}\mathbf{l}$, $[{\mathbf{F_p}}(\mathbf{jl})]^*=-{\mathbf{F_p}}(\mathbf{jl})$ and $[{\mathbf{F_p}}(\mathbf{l}+\mathbf{jl})/\sqrt{2}]^*=-{\mathbf{F_p}}(\mathbf{l}+\mathbf{jl})/\sqrt{2}$. By the multiplication rule shown in \eqref{eq:relation-remark}, it is easy to derive $$\mathbf{OO^*}=-\mathbf{OO}=-({\mathbf{F_p}\mathbf{p}})({\mathbf{F_p}\mathbf{p}})
={\mathbf{F_p}}{\mathbf{F_p^*}}=\mathbf{I_p}.$$
Hence, octonion matrices ${\mathbf{F_p}}\mathbf{l}$, ${\mathbf{F_p}}(\mathbf{jl})$ and ${\mathbf{F_p}}(\mathbf{l}+\mathbf{jl})/\sqrt{2}$ are unitary. \ep

Second, we prove that unitary octonion matrices ${\mathbf{F_p}}\mathbf{l}$, ${\mathbf{F_p}}(\mathbf{jl})$ and ${\mathbf{F_p}}(\mathbf{l}+\mathbf{jl})/\sqrt{2}$ are all diagonalization matrices of any given $p\times p$ circulant quaternion matrix.
\begin{theorem}\label{theo:3}
Let $\mathbf{p}=\mathbf{l}, \mathbf{jl}$ or $ (\mathbf{l}+\mathbf{jl})/\sqrt{2}$, and $\mathbf{O}={\mathbf{F_p}}\mathbf{p}$. Then $\mathbf{O}$ is a diagonalization matrices of $p\times p$ circulant quaternion matrix, i.e.,
for any circulant quaternion matrix $\mathbf{{Q}}=\operatorname{circ}(q_1,q_2,\ldots,q_p)\in \mathbb{Q}^{p\times p}$,
$${\mathbf{O}}\mathbf{{Q}}{\mathbf{O}}^*={\rm Diag}(\hat{q}_1,\hat{q}_2,\ldots,\hat{q}_p)\in \mathbb{Q}^{p\times p}.$$
Further, let $\hat{\mathbf{q}}:=(\hat{q}_1,\hat{q}_2,\ldots,\hat{q}_p)^\top$ and $\mathbf{q}:=(q_1,q_2,\ldots,q_p)^\top$, then $\hat{\mathbf{q}}$ and $\mathbf{q}$ satisfy:
\begin{equation}\label{relation--1}
\hat{\mathbf{q}}=\sqrt{p}\mathbf{F_{p}}\bar{\mathbf{q}}.
\end{equation}
%
\end{theorem}

{\bf Proof.} Rewrite the circulant quaternion matrix $\mathbf{{Q}}=\operatorname{circ}(q_1,q_2,\ldots,q_p)$ by the CD from  $\mathbf{Q}=\mathbf{Q_1}+\mathbf{Q_2}\mathbf{j}$, where $\mathbf{Q_1}$ and $ \mathbf{Q_2}$ are circulant complex matrices. Denote
$$
\mathbf{Q_1}= {\rm circ}(q_1^1,q_2^1,\ldots,q_p^1),\  \mathbf{Q_2}={\rm circ}(q_1^2,q_2^2,\ldots,q_p^2),
$$
where $q_s^j\in \mathbb{C}$ for every $s\in [p]$ and $j=1,2$. Then
$ q_s= q_s^1 + q_s^2 \mathbf{j} $.

By the multiplication rule shown in \eqref{eq:relation-remark}, it is deduced
\begin{equation}\label{eq:O=O=O}
{\mathbf{O}}\mathbf{{Q}}{\mathbf{O}}^*
=-({\mathbf{F_p}\mathbf{p}})\mathbf{Q}({\mathbf{F_p}\mathbf{p}})
=-({\mathbf{F_p}\mathbf{p}})(\mathbf{Q_1}+\mathbf{Q_2}\mathbf{j})({\mathbf{F_p}\mathbf{p}})
=\mathbf{F_p}\overline{\mathbf{Q_1}}\mathbf{F_p}^*-
\mathbf{F_p}{\mathbf{Q_2}}\mathbf{F_p}^*\mathbf{j}.
\end{equation}
Since $\mathbf{Q_1}$ and $ \mathbf{Q_2}$ are circulant complex matrices and $\mathbf{F_p}$ is the diagonalization matrix of $p\times p$ circulant complex matrix, it follows that $\mathbf{F_p}\overline{\mathbf{Q_1}}\mathbf{F_p^*}$ and $\mathbf{F_p}{\mathbf{Q_2}}\mathbf{F_p^*}$ are diagonal complex matrices. Thus, we have that ${\mathbf{O}}\mathbf{{Q}}{\mathbf{O}}^*$ is diagonal, denote by ${\rm Diag}(\hat{q}_1,\hat{q}_2,\ldots,\hat{q}_p).$
Denote
\begin{equation}\label{eq:DiagFQF}
{\rm Diag}(\hat{q}_1^1,\hat{q}_2^1,\ldots,\hat{q}_p^1)=
\mathbf{F_p}\overline{\mathbf{Q_1}}\mathbf{F_p^*},\quad\quad{\rm Diag}(\hat{q}_1^2,\hat{q}_2^2,\ldots,\hat{q}_p^2)=
\mathbf{F_p}{\mathbf{Q_2}}\mathbf{F_p^*}.
\end{equation}
Then it follows that
\begin{equation}\label{eq:hatq}
\hat{q}_s^1=\sum\limits_{m=1}^p\omega^{(s-1)(m-1)}\overline{q_{m}^1}\quad{\rm and}\quad \hat{q}_s^2=\sum\limits_{m=1}^p\omega^{(s-1)(m-1)}q_{m}^2,
\end{equation}
for all $s\in [p]$. Further, from \eqref{eq:O=O=O}, \eqref{eq:DiagFQF} and \eqref{eq:hatq}, we have that for all $s\in [p]$,
\begin{equation}\label{eq:hatq}
\hat{q}_s=\hat{q}_s^1-\hat{q}_s^2\mathbf{j}=\sum\limits_{m=1}^p\omega^{(s-1)(m-1)}\overline{q_{m}},
\end{equation}
which implies that \eqref{relation--1} holds. This completes the proof. \ep

\begin{remark}\label{rem-1}
\mbox{\rm (i)} Recall that if two complex vector $\mathbf{{a}},\mathbf{{b}}\in \mathbb{C}^p$ satisfy that $\mathbf{{a}}={\mathbf{F_p}}\mathbf{{b}}$, then $\mathbf{{a}}$ can be obtained by applying the FFT to $\mathbf{{b}}$, denoted by $\mathbf{{a}}=\operatorname{fft}(\mathbf{{b}})$. Therefore, we can compute $\hat{\mathbf{{q}}}$ in \eqref{relation--1} by applying the FFT to $\bar{\mathbf{{q}}}=\overline{\mathbf{q_1}}-\mathbf{q_2}\mathbf{j}$ with $\mathbf{q_1},\mathbf{q_2}\in \mathbb{C}^p$, i.e., \begin{equation}\label{eq:fft-q}
\hat{\mathbf{{q}}}=\operatorname{fft}(\overline{\mathbf{q_1}})-
\operatorname{fft}(\mathbf{{q_2}})\mathbf{j},
\end{equation}
denoted by $\hat{\mathbf{{q}}}=\operatorname{fft}(\bar{\mathbf{{q}}}).$

\mbox{\rm (ii)} For any ${q}\in \mathbb{C}^p$, the cost of the operation $\operatorname{fft}({q})$ is $O(p\log p)$ \cite{GV-2013}, thus from \eqref{eq:fft-q}, the cost of diagonalization of any ${p\times p}$ circulant quaternion matrix is $O(p\log p)$ by the way given in Theorem \ref{theo:3}.

\mbox{\rm (iii)} Let $\mathbf{p}=\mathbf{l}, \mathbf{jl}$ or $ (\mathbf{l}+\mathbf{jl})/\sqrt{2}$, and $\mathbf{O}={\mathbf{F^*_p}}\mathbf{p}$. By a similar process as in Theorem \ref{theo:3}, it can be proved that for any $\hat{\mathbf{{q}}}=(\hat{q}_1,\hat{q}_2,\ldots,\hat{q}_p)^\top\in \mathbb{Q}^p$ and $\mathbf{{q}}=(q_1,q_2,\ldots,q_p)^\top\in \mathbb{Q}^p$, if
$$\operatorname{circ}(q_1,q_2,\ldots,q_p)={\mathbf{O}}{\rm Diag}(\hat{q}_1,\hat{q}_2,\ldots,\hat{q}_p){\mathbf{O}^*},$$
then $\hat{\mathbf{{q}}}$ and $\mathbf{q}$ satisfy:
$\mathbf{{q}}=\frac{1}{\sqrt{p}}\mathbf{F^*_{p}}\bar{\hat{\mathbf{q}}},$ which means $\mathbf{{q}}$ can be computed by applying the inverse FFT (ifft) to $\bar{\hat{\mathbf{q}}}=\overline{{\mathbf{\hat{q}_1}}}- {\mathbf{\hat{q}_2}} \mathbf{j}$ with $\mathbf{\hat{q}_1},\mathbf{\hat{q}_2}\in \mathbb{C}^p$, i.e.,
$$
\mathbf{{q}}=\operatorname{ifft}(\overline{\mathbf{\hat{q}_1}})-
\operatorname{ifft}({\mathbf{\hat{q}_2}})\mathbf{j},
$$
denoted by ${\mathbf{q}}=\operatorname{ifft}({\bar{\hat{\mathbf{q}}}})$.
\end{remark}

\subsection{Block diagonalization of block circulant quaternion matrices}
\label{subsect.blockdiag-oct}
It has been shown in Subsection \ref{Sect.2} that for any block circulant matrix $\mathbf{{A}}={\rm bcirc}(\mathbf{A_1},\mathbf{A_2},\ldots,\mathbf{A_p})\in \mathbb{C}^{mp\times np}$, it can be block diagonalized by unitary matrices $\mathbf{F_p}\otimes \mathbf{I_m}\in \mathbb{C}^{mp\times mp}$ and $\mathbf{F_p}\otimes \mathbf{I_n} \in \mathbb{C}^{np\times np}$
as
\begin{equation}\label{eq:complex-block-diag}
(\mathbf{F_p}\otimes \mathbf{I_m}) \mathbf{{A}} (\mathbf{F^*_p}\otimes \mathbf{I_n}) ={\rm Diag}(\mathbf{\hat{A}_1},\mathbf{\hat{A}_2},\ldots,\mathbf{\hat{A}_p})\in \mathbb{C}^{mp\times np},
\end{equation}
where for all $s\in[p]$,
\begin{equation}\label{eq:complex-relation-block}
\mathbf{\hat{A}_s}=\sum\limits_{r=1}^p\omega^{(s-1)(r-1)}\mathbf{A_r}.
\end{equation}
Here, we prove that any block circulant quaternion matrix $\mathbf{{Q}}\in \mathbb{Q}^{mp\times np}$ can be block diagonalized by unitary matrices $\mathbf{F_p}\mathbf{p}\otimes \mathbf{I_m}\in \mathbb{O}^{mp\times mp}$ and $\mathbf{F_p}\mathbf{p}\otimes \mathbf{I_n} \in \mathbb{O}^{np\times np}$, where $\mathbf{p}=\mathbf{l}, \mathbf{jl}$ or $ (\mathbf{l}+\mathbf{jl})/\sqrt{2}$.

\begin{theorem}\label{theo:block-diagonal}
Let $\mathbf{p}=\mathbf{l}, \mathbf{jl}$ or $ (\mathbf{l}+\mathbf{jl})/\sqrt{2}$, and $\mathbf{O}={\mathbf{F_p}}\mathbf{p}$.
Then ${\mathbf{O}}\otimes \mathbf{I_m}$ is a unitary octonion matrix and for any block circulant quaternion matrix $\mathbf{{Q}}=\operatorname{bcirc}(\mathbf{{Q_1}},\mathbf{{Q_2}},\ldots,\mathbf{{Q_p}})\in \mathbb{Q}^{mp\times np}$ with each $\mathbf{{Q}_r}\in \mathbb{Q}^{m\times n}$,
\begin{equation}\label{eq:block-diagonal}
({\mathbf{O}}\otimes \mathbf{I_m})\mathbf{{Q}}({\mathbf{O}^*}\otimes \mathbf{I_n}) = {\rm Diag}({\mathbf{\hat{Q}_1}},{\mathbf{\hat{Q}_2}},\ldots,{\mathbf{\hat{Q}_p}})\in \mathbb{Q}^{mp\times np},
\end{equation}
where each
\begin{equation}\label{relation--2}
\mathbf{\hat{Q}_s}=\sum\limits_{r=1}^p\omega^{(s-1)(r-1)}\overline{\mathbf{Q_r}}\in \mathbb{Q}^{m\times n}.
\end{equation}
\end{theorem}

{\bf Proof.} First, by the definion of the Kronecker product ``$\otimes$'' in \eqref{kronecker}, we have that
$$({\mathbf{O}}\otimes \mathbf{I_m})^*={\mathbf{O}}^*\otimes \mathbf{I_m}.$$
Further, from Theorem \ref{ther:unitary} and the property of that ``$(\mathbf{A}\otimes \mathbf{B})(\mathbf{C}\otimes \mathbf{D})=\mathbf{AC}\otimes \mathbf{BD}$'', it follows that
\begin{eqnarray*}
 ({\mathbf{O}}\otimes \mathbf{I_m})({\mathbf{O}}\otimes \mathbf{I_m})^* = ({\mathbf{O}}\mathbf{O}^*)\otimes (\mathbf{I_m}\mathbf{I_m})=\mathbf{I_{mp}},
\end{eqnarray*}
which means that ${\mathbf{O}}\otimes \mathbf{I_m}$ is unitary.

Second, we prove that \eqref{eq:block-diagonal} and \eqref{relation--2} hold. Denote the entry in the $i$-th row, $j$-th column of $\mathbf{{Q_r}}$ as $q^r_{ij}$ for all $i\in [m]$, $j\in [n]$ and $r\in[p]$, then it follows that
$$\mathbf{{Q}}=\operatorname{bcirc}(\mathbf{{Q_1}},\mathbf{{Q_2}},\ldots,\mathbf{{Q_p}})
=\sum_{i=1}^m\sum_{j=1}^n[{\rm circ}(q^1_{ij},q^2_{ij},\ldots,q^p_{ij})]\otimes \mathbf{E_{ij}},$$
where $\mathbf{E_{ij}}$ means the matrix whose entries are zeros except the $(i,j)$-th entry being one. Thus, it is deduced that
\begin{eqnarray}
\nonumber({\mathbf{O}}\otimes \mathbf{I_m})\mathbf{{Q}}({\mathbf{O}^*}\otimes \mathbf{I_n})&=&({\mathbf{O}}\otimes \mathbf{I_m})\left(\sum_{i=1}^m\sum_{j=1}^n[{\rm circ}(q^1_{ij},q^2_{ij},\ldots,q^p_{ij})]\otimes \mathbf{E_{ij}}\right)({\mathbf{O}^*}\otimes \mathbf{I_n})\\
\nonumber&=&\sum_{i=1}^m\sum_{j=1}^n \left({\mathbf{O}}[{\rm circ}(q^1_{ij},q^2_{ij},\ldots,q^p_{ij}]{\mathbf{O}^*}\right)\otimes (\mathbf{I_m}\mathbf{E_{ij}}\mathbf{I_n})\\
&=&\sum_{i=1}^m\sum_{j=1}^n \left({\mathbf{O}}[{\rm circ}(q^1_{ij},q^2_{ij},\ldots,q^p_{ij}]{\mathbf{O}^*}\right)\otimes \mathbf{E_{ij}}.\label{eq:kron-big}
\end{eqnarray}

In Theorem \ref{theo:3}, it is shown that each circulant quaternion matrix ${\rm circ}(q^1_{ij},q^2_{ij},\ldots,q^p_{ij})$ can be diagonalized by $\mathbf{O}$ and $\mathbf{O}^*$. That is to say, each ${\mathbf{O}}[{\rm circ}(q^1_{ij},q^2_{ij},\ldots,q^p_{ij}]{\mathbf{O}^*}$ is a diagonal quaternion matrix. Thus, by the Kronecker product ``$\otimes$'' shown in \eqref{kronecker}, it follows
$({\mathbf{O}}\otimes \mathbf{I_m})\mathbf{{Q}}({\mathbf{O}^*}\otimes \mathbf{I_n})=\sum\limits_{i=1}^m\sum\limits_{j=1}^n \left({\mathbf{O}}[{\rm circ}(q^1_{ij},q^2_{ij},\ldots,q^p_{ij}]{\mathbf{O}^*}\right)\otimes \mathbf{E}_{ij}$ is a block diagonal quatenion matrix.

Furthermore, let $({\mathbf{O}}\otimes \mathbf{I_m})\mathbf{{Q}}({\mathbf{O}^*}\otimes \mathbf{I_n}) = {\rm Diag}({\mathbf{\hat{Q}_1}},\mathbf{\hat{Q}_2},\ldots, \mathbf{\hat{Q}_p})$ as defined in \eqref{eq:block-diagonal}, and denote the entry in the $i$-th row, $j$-th column of ${\mathbf{\hat{Q}_r}}$ as $\hat{q}^r_{ij}$ for all $i\in [m]$, $j\in [n]$ and $r\in[p]$. Then it follows that
$$[\mathbf{\hat{Q}_1^\top},\mathbf{\hat{Q}_2^\top},\ldots,\mathbf{\hat{Q}_p^\top}]^\top
=\sum_{i=1}^m\sum_{j=1}^n
[\hat{q}^1_{ij},\hat{q}^2_{ij},\ldots,\hat{q}^p_{ij}]^\top\otimes \mathbf{E_{ij}},$$
which implies that these $\hat{q}^1_{ij},\hat{q}^2_{ij},\ldots,\hat{q}^p_{ij}$ are exact the diagonal entries of ${\mathbf{O}}[{\rm circ}(q^1_{ij},q^2_{ij},\ldots,q^p_{ij}]{\mathbf{O}^*}$ by \eqref{eq:kron-big}.
From Theorem \ref{theo:3}, we can obtain that for any $i\in[m]$, $j\in [n]$ and $s\in [p]$, each
$$\hat{q}^s_{ij}=\sum\limits_{r=1}^p\omega^{(s-1)(r-1)}\overline{{q}^r_{ij}}.$$
Noting that for every $s\in [p]$,
$$\mathbf{\hat{Q}_s}=\sum_{i=1}^m\sum_{j=1}^n\hat{q}^s_{ij}\otimes \mathbf{E_{ij}},\quad{\rm and}\quad \mathbf{{Q}_s}=\sum_{i=1}^m\sum_{j=1}^n{q}^s_{ij}\otimes \mathbf{E_{ij}},$$
thus, we obtain that
$$\mathbf{\hat{Q}_s}=\sum\limits_{r=1}^p\omega^{(s-1)(r-1)}\overline{\mathbf{Q_r}},$$
i.e., \eqref{relation--2} holds. This completes the proof.\ep

Recall that it is shown in Subsection \ref{Sub:T-product} that block circulant matrices have a one-by-one correspondence with third-order tensors, by the operator ``$bcirc$''. Hence, we can immediately obtain the following tensor version of Theorem \ref{theo:block-diagonal}.

\begin{theorem}\label{Coro:tensor}
Let $\mathcal{Q}=\mathcal{Q}_1+\mathcal{Q}_2\mathbf{j}\in \mathbb{Q}^{m\times n\times p}$ with $\mathcal{Q}_1,\mathcal{Q}_2\in \mathbb{C}^{m\times n\times p}$, $\mathbf{p}=\mathbf{l}, \mathbf{jl}$ or $ (\mathbf{l}+\mathbf{jl})/\sqrt{2}$, and $\mathbf{O}={\mathbf{F_p}}\mathbf{p}$. Then the block circulant matrix $bcirc(\mathcal{Q})\in \mathbb{Q}^{mp\times np}$
can be block diagonalized by $\mathbf{O}\otimes \mathbf{I_m }\in \mathbb{O}^{mp\times mp}$ and $\mathbf{O}^*\otimes \mathbf{I_n} \in \mathbb{O}^{np\times np}$, that is
\begin{equation}\label{eq:block-tensor}
(\mathbf{O}\otimes \mathbf{I_m})\ bcirc(\mathcal{Q})\ (\mathbf{O}^*\otimes \mathbf{I_n})={\rm Diag}(\hat{\mathcal{Q}})\in \mathbb{Q}^{mp\times np}
\end{equation}
where $\hat{\mathcal{Q}}$  and $\mathcal{Q}$ satisfy:
\begin{equation}\label{eq:relation-unfold}
unfold(\hat{\mathcal{Q}})=\sqrt{p}(\mathbf{F_p}\otimes \mathbf{I_m})unfold(\overline{\mathcal{Q}}).
\end{equation}
That means, the quaternion tensor $\hat{\mathcal{Q}}$ can be computed by applying the fast fourier transform to $\overline{\mathcal{Q}}=\overline{\mathcal{Q}_1}-\mathcal{Q}_2\mathbf{j}$ along the third dimension as:
\begin{equation}\label{eq:fft-Q}
\hat{\mathcal{Q}}=\operatorname{fft}(\overline{\mathcal{Q}},[\;],3):=\operatorname{fft}(\overline{\mathcal{Q}_1},[\;],3) -\operatorname{fft}({\mathcal{Q}}_2,[\;],3)\mathbf{j}.
\end{equation}
\end{theorem}

{\bf Proof.} From \eqref{eq:block-diagonal} and \eqref{relation--2} in Theorem \ref{theo:block-diagonal}, it follows that \eqref{eq:block-tensor} and \eqref{eq:relation-unfold} hold immediately. By the Kronecker product ``$\otimes$'' defined in \eqref{kronecker} and the ``$unfold$'' operator defined in \eqref{bcirc-unfold}, we can further obtain that for any $i\in [m]$ and $j\in [n]$,
$$[(\mathbf{\hat{Q}^{(1)}})_{ij},(\mathbf{\hat{Q}^{(2)}})_{ij},\ldots,
(\mathbf{\hat{Q}^{(p)}})_{ij}]^\top
=\sqrt{p}\mathbf{F_{p}}[(\overline{\mathbf{Q^{(1)}}})_{ij},(\overline{\mathbf{Q^{(2)}}})_{ij},
\ldots,(\overline{\mathbf{Q^{(p)}}})_{ij}]^\top.$$
This implies that the quaternion tensor $\hat{\mathcal{Q}}$ can be computed by \eqref{eq:fft-Q}. \ep

\begin{remark}\label{rem-2}
\mbox{\rm (i)} From Theorems \ref{theo:3}, \ref{theo:block-diagonal} and  \ref{Coro:tensor}, it can be seen that although the diagonalization of the circulant quaternion matrix or the block diagonalization of the block circulant quaternion matrix is achieved by unitary octonion matrices, the result diagonal matrix or block diagonal matrix remains in the quaternion domain.

%
\mbox{\rm (ii)} The cost of the method that achieves the block diagonalization of any ${mp\times np}$ block circulant quaternion matrix $\mathbf{Q}$ via $\operatorname{fft}(\overline{\mathcal{Q}},[\;],3)$ as defined by \eqref{eq:fft-Q} is $O(mnp\log p)$, since for any $\mathcal{Q}\in \mathbb{C}^{m\times n\times p}$, the cost of the operation $\operatorname{fft}(\mathcal{Q},[\;],3)$ is $O(mnp\log p)$ \cite{KM-2011}.

\mbox{\rm (iii)} Let $\mathbf{p}=\mathbf{l}, \mathbf{jl}$ or $ (\mathbf{l}+\mathbf{jl})/\sqrt{2}$. Suppose that $\hat{\mathcal{Q}}$, $\mathcal{Q}\in\mathbb{Q}^{m\times n\times p}$ satisfy
$$(\mathbf{F_p}\mathbf{p}\otimes \mathbf{I_m}) bcirc(\mathcal{Q}) (-\mathbf{F_p}\mathbf{p}\otimes \mathbf{I_n})={\rm Diag}(\hat{\mathcal{Q}}).$$ It should be noted that different from the case of $\mathcal{Q}\in \mathbb{C}^{m\times n\times p}$, it is deduced that
$${\rm Diag}(\hat{\mathcal{Q}})\neq(\mathbf{F}^*_{p}\mathbf{p}\otimes \mathbf{I}_{m})bcirc(\mathcal{Q})(-\mathbf{F}^*_{p}\mathbf{p}\otimes\mathbf{I}_{n})$$
in general, due to the non-associativity of octonion algebra.
Actually, from Proposition \ref{pro:associative} {\rm (iii)} and {\rm (iv)}, it follows that
\begin{equation}\label{eq:ifftA}
\begin{array}{rcl}
&&bcirc(\mathcal{Q})\\
&=&(-\mathbf{F_{p}}\mathbf{p}\otimes \mathbf{I_{m}}){\rm Diag}(\hat{\mathcal{Q}}_1)(\mathbf{F_{p}}\mathbf{p}\otimes\mathbf{I_{n}})
+(-\mathbf{F^*_{p}}\mathbf{p}\otimes \mathbf{I_{m}})[{\rm Diag}(\hat{\mathcal{Q}}_2)\mathbf{j}](\mathbf{F^*_{p}}\mathbf{p}\otimes\mathbf{I_{n}})\\
&=&(-\mathbf{F^*_{p}}\mathbf{p}\otimes \mathbf{I_{m}})[(\mathbf{{P}_{p}}\otimes \mathbf{I_{m}}){\rm Diag}({\hat{\mathcal{Q}}_1})(\mathbf{{P}_{p}}\otimes\mathbf{I_{n}})+{\rm Diag}(\hat{\mathcal{Q}}_2)\mathbf{j}](\mathbf{F^*_{p}}\mathbf{p}\otimes\mathbf{I_{n}}),
\end{array}
\end{equation}
where $\mathbf{{P}_{p}}=\mathbf{F_{p}}\mathbf{F_{p}}=\mathbf{F^*_{p}}\mathbf{F^*_{p}} \in \mathbb{R}^{p\times p}$ is a permutation matrix defined by \eqref{eq:Pp}.

\mbox{\rm (iv)} Let $\mathbf{p}=\mathbf{l}, \mathbf{jl}$ or $ (\mathbf{l}+\mathbf{jl})/\sqrt{2}$, and $\mathbf{O}={\mathbf{F^*_p}}\mathbf{p}$. With Remark \ref{rem-1}, by a similar process as in Theorems \ref{theo:block-diagonal} and  \ref{Coro:tensor}, it can be also proved that for third-order tensors $\hat{\mathcal{Q}}$ and $\mathcal{Q}$, if
$$bcirc(\mathcal{Q})=(\mathbf{O}\otimes \mathbf{I_m})\ {\rm Diag}(\hat{\mathcal{Q}})\ (\mathbf{O}^*\otimes \mathbf{I_n}),$$
then $\hat{\mathcal{Q}}$  and $\mathcal{Q}$ satisfy:
$$
unfold({\mathcal{Q}})=\frac{1}{\sqrt{p}}(\mathbf{F^*_p}\otimes \mathbf{I_m})unfold(\overline{\hat{\mathcal{Q}}}),
$$
which means $\mathcal{Q}$ can be computed by applying the inverse fast fourier transform to $\overline{\hat{\mathcal{Q}}}=\overline{\mathcal{Q}_1}-\mathcal{Q}_2\mathbf{j}$ along the third dimension as:
$$
\mathcal{Q}=\operatorname{ifft}(\overline{\hat{\mathcal{Q}}},[\;],3)
:=\operatorname{ifft}(\overline{\hat{\mathcal{Q}}_1},[\;],3) -\operatorname{ifft}(\hat{\mathcal{Q}}_2,[\;],3)\mathbf{j}.
$$
\end{remark}

\section{The T-product of third-order quaternion tensors}
\label{Sect.4}
\setcounter{equation}{0} \setcounter{assumption}{0}
\setcounter{theorem}{0} \setcounter{proposition}{0}
\setcounter{corollary}{0} \setcounter{lemma}{0}
\setcounter{definition}{0} \setcounter{remark}{0}
\setcounter{algorithm}{0}
In this section, we apply the block diagonalization of block circulant quaternion matrices achieved in Subsection \ref{subsect.blockdiag-oct} to calculate the T-product between quaternion tensors and give a positive answer to {\bf Question 3}.  By Theorem \ref{Coro:tensor}, for every third-order quaternion tensor $\mathcal{Q}\in \mathbb{Q}^{m\times n\times p}$,
$$
(\mathbf{O}\otimes \mathbf{I_m})\ bcirc(\mathcal{Q})\ (\mathbf{O}^*\otimes \mathbf{I_n})={\rm Diag}(\hat{\mathcal{Q}})\in \mathbb{Q}^{mp\times np}$$
where $\mathbf{O}=\mathbf{F_p}\mathbf{l}$, $\mathbf{F_p}(\mathbf{jl})$ or  $\mathbf{F_p}(\mathbf{l}+\mathbf{jl})/\sqrt{2}$. Here, we choose the case of $\mathbf{O}=\mathbf{F_p}\mathbf{l}$ to finish the following discussion.
The complex form of the quaternion tensor $\mathcal{A}$ (i.e., $\mathcal{A}=\mathcal{A}_1+\mathcal{A}_2\mathbf{j}$ with $\mathcal{A}_1,\mathcal{A}_2$ being complex) plays an important role in the following discussions. In the rest of paper, for simplicity, we use the notation $\mathcal{A}_1$ and $\mathcal{A}_2$ to represent the complex tensors $\mathcal{A}_1$ and $\mathcal{A}_2$ in the complex form of the quaternion tensor $\mathcal{A}$, unless otherwise specified.

\begin{theorem}\label{theo:T-product}
Let $\mathcal{A}\in \mathbb{Q}^{m\times n\times p}$, $\mathcal{B}\in \mathbb{Q}^{n\times s\times p}$, $\operatorname{fft}(\overline{\mathcal{A}},[\;],3)$ and $\operatorname{fft}(\overline{\mathcal{B}},[\;],3)$ be defined by \eqref{eq:fft-Q}. Denote
\begin{equation}\label{eq:fftA-B}
\hat{\mathcal{A}}=\hat{\mathcal{A}}_1+\hat{\mathcal{A}}_2\mathbf{j} =\operatorname{fft}(\overline{\mathcal{A}},[\;],3),\quad{\rm and}\quad\hat{\mathcal{B}}=\hat{\mathcal{B}}_1+\hat{\mathcal{B}}_2\mathbf{j} =\operatorname{fft}(\overline{\mathcal{B}},[\;],3).
\end{equation}
Then the T-product $\mathcal{C}=\mathcal{A}\ast\mathcal{B}$ can be obtained by
\begin{equation}\label{eq:T-product-quaternion}
\mathcal{C}=bcirc^{-1}[(\mathbf{F^*_{p}}\mathbf{l}\otimes \mathbf{I_{m}}){\rm Diag}({\hat{\mathcal{C}}})(-\mathbf{F^*_{p}}\mathbf{l}\otimes\mathbf{I_{s}})],
\end{equation}
where ${\rm Diag}({\hat{\mathcal{C}}})={\rm Diag}({\hat{\mathcal{C}}_1})+{\rm Diag}({\hat{\mathcal{C}}_2})\mathbf{j}$,
\begin{eqnarray}
{\rm Diag}(\hat{\mathcal{C}}_1)&=&(\mathbf{{P}_{p}}\otimes \mathbf{I_{m}}){\rm Diag}({\hat{\mathcal{A}}_1}){\rm Diag}({\hat{\mathcal{B}}_1})(\mathbf{{P}_{p}}\otimes\mathbf{I_{s}})\label{C}\\
&&-{\rm Diag}(\overline{{\hat{\mathcal{A}}_2}})[(\mathbf{{P}_{p}}\otimes \mathbf{I_{n}}){\rm Diag}(\hat{\mathcal{B}}_2)(\mathbf{{P}_{p}}\otimes\mathbf{I_{s}})],\nonumber \\
{\rm Diag}({\hat{\mathcal{C}}_2})&=&[(\mathbf{{P}_{p}}\otimes \mathbf{I_{m}}){\rm Diag}(\overline{\hat{\mathcal{A}}_1})(\mathbf{{P}_{p}}\otimes\mathbf{I_{n}})]{\rm Diag}(\hat{\mathcal{B}}_2)\label{C2}\\
&&+{\rm Diag}(\hat{\mathcal{A}}_2){\rm Diag}(\hat{\mathcal{B}}_1),\nonumber
\end{eqnarray}
and $\mathbf{{P}_{p}}\in \mathbb{R}^{p\times p}$ is a permutation matrix defined by \eqref{eq:Pp}.
\end{theorem}

{\bf Proof.} To show that \eqref{eq:T-product-quaternion} holds, we only need to prove that
\begin{equation}\label{eq:T-product}
bcirc(\mathcal{C})=(\mathbf{F^*_{p}}\mathbf{l}\otimes \mathbf{I_{m}})[{\rm Diag}({\hat{\mathcal{C}}_1})+{\rm Diag}({\hat{\mathcal{C}}_2})\mathbf{j}](-\mathbf{F^*_{p}}\mathbf{l}\otimes\mathbf{I_{s}}).
\end{equation}
From Theorem \ref{Coro:tensor} and \eqref{eq:fftA-B}, it follows that
\begin{equation}\label{eq:A}
{\rm Diag}(\hat{\mathcal{A}})={\rm Diag}(\hat{\mathcal{A}}_1)+{\rm Diag}(\hat{\mathcal{A}}_2)\mathbf{j}=(\mathbf{F_{p}}\mathbf{l}\otimes \mathbf{I_{m}})bcirc(\mathcal{A})(-\mathbf{F_{p}}\mathbf{l}\otimes\mathbf{I_{n}}),
\end{equation}
\begin{equation}\label{eq:B}
{\rm Diag}(\hat{\mathcal{B}})={\rm Diag}(\hat{\mathcal{B}}_1)+{\rm Diag}(\hat{\mathcal{B}}_2)\mathbf{j}=(\mathbf{F_{p}}\mathbf{l}\otimes \mathbf{I_{m}})bcirc(\mathcal{B})(-\mathbf{F_{p}}\mathbf{l}\otimes\mathbf{I_{n}}).
\end{equation}
Then by Proposition \ref{pro:associative} (iii), (iv) and \eqref{eq:A}, it follows that
\begin{equation}\label{eq:AAA}
bcirc(\mathcal{A})=(-\mathbf{F_{p}}\mathbf{l}\otimes \mathbf{I_{m}}){\rm Diag}(\hat{\mathcal{A}}_1)(\mathbf{F_{p}}\mathbf{l}\otimes\mathbf{I_{n}})+(-\mathbf{F^*_{p}}\mathbf{l}\otimes \mathbf{I_{m}})[{\rm Diag}(\hat{\mathcal{A}}_2)\mathbf{j}](\mathbf{F^*_{p}}\mathbf{l}\otimes\mathbf{I_{n}}).
\end{equation}
By Proposition \ref{pro:associative} (iii), (iv) and \eqref{eq:B}, it follows that
\begin{equation}\label{eq:BBB}bcirc(\mathcal{B})=(-\mathbf{F_{p}}\mathbf{l}\otimes \mathbf{I_{n}}){\rm Diag}(\hat{\mathcal{B}}_1)(\mathbf{F_{p}}\mathbf{l}\otimes\mathbf{I_{s}})+(-\mathbf{F^*_{p}}\mathbf{l}\otimes \mathbf{I_{n}})[{\rm Diag}(\hat{\mathcal{B}}_2)\mathbf{j}](\mathbf{F^*_{p}}\mathbf{l}\otimes\mathbf{I_{s}}).
\end{equation}
By a similar proof as \cite[Lemma 2.3]{L-2020}, we can obtain that the following property holds for quaternion tensors $\mathcal{A}$ and $\mathcal{B}$, that is:
\begin{equation}\label{eq:bcirc}
bcirc(\mathcal{A}\ast\mathcal{B})=bcirc(\mathcal{A})\cdot bcirc(\mathcal{B}).
\end{equation}
Thereby, from \eqref{eq:AAA}, \eqref{eq:BBB} and \eqref{eq:bcirc}, it is deduced that
\begin{equation}\label{bcirc(C)}
\begin{array}{rl}
&bcirc(\mathcal{C})=bcirc(\mathcal{A}\ast\mathcal{B})=bcirc(\mathcal{A})\cdot bcirc(\mathcal{B})\\
=&[(-\mathbf{F_{p}}\mathbf{l}\otimes \mathbf{I_{m}}){\rm Diag}(\hat{\mathcal{A}}_1)(\mathbf{F_{p}}\mathbf{l}\otimes\mathbf{I_{n}})]\cdot [(-\mathbf{F_{p}}\mathbf{l}\otimes \mathbf{I_{n}}){\rm Diag}(\hat{\mathcal{B}}_1)(\mathbf{F_{p}}\mathbf{l}\otimes\mathbf{I_{s}})]\\
&+\;[(-\mathbf{F^*_{p}}\mathbf{l}\otimes \mathbf{I_{m}})({\rm Diag}(\hat{\mathcal{A}}_2)\mathbf{j})(\mathbf{F^*_{p}}\mathbf{l}\otimes\mathbf{I_{n}})]\cdot
[(-\mathbf{F^*_{p}}\mathbf{l}\otimes \mathbf{I_{n}})({\rm Diag}(\hat{\mathcal{B}}_2)\mathbf{j})(\mathbf{F^*_{p}}\mathbf{l}\otimes\mathbf{I_{s}})]\\
&+\;[(-\mathbf{F_{p}}\mathbf{l}\otimes \mathbf{I_{m}}){\rm Diag}(\hat{\mathcal{A}}_1)(\mathbf{F_{p}}\mathbf{l}\otimes\mathbf{I_{n}})]\cdot[(-\mathbf{F^*_{p}}\mathbf{l}\otimes \mathbf{I_{n}})({\rm Diag}(\hat{\mathcal{B}}_2)\mathbf{j})(\mathbf{F^*_{p}}\mathbf{l}\otimes\mathbf{I_{s}})]\\
&+\;[(-\mathbf{F^*_{p}}\mathbf{l}\otimes \mathbf{I_{m}})({\rm Diag}(\hat{\mathcal{A}}_2)\mathbf{j})(\mathbf{F^*_{p}}\mathbf{l}\otimes\mathbf{I_{n}})]\cdot
[(-\mathbf{F_{p}}\mathbf{l}\otimes \mathbf{I_{n}}){\rm Diag}(\hat{\mathcal{B}}_1)(\mathbf{F_{p}}\mathbf{l}\otimes\mathbf{I_{s}})].
\end{array}
\end{equation}
By \eqref{eq:relation-remark} and \eqref{FpFp}, it is easy to see $\mathbf{{P}_{p}}=\mathbf{F_{p}}\mathbf{F_{p}}=\mathbf{F^*_{p}}\mathbf{F^*_{p}} =-(\mathbf{F_{p}}\mathbf{l})(\mathbf{F^*_{p}}\mathbf{l})
=-(\mathbf{F^*_{p}}\mathbf{l})(\mathbf{F_{p}}\mathbf{l})\in \mathbb{R}^{p\times p}$ and
\begin{eqnarray}
\mathbf{{P}_{p}}\otimes \mathbf{I_m}
=-(\mathbf{F_{p}}\mathbf{l}\otimes \mathbf{I_m})(\mathbf{F^*_{p}}{\mathbf{l}}\otimes \mathbf{I_m})=-(\mathbf{F^*_{p}}\mathbf{l}\otimes \mathbf{I_m})(\mathbf{F_{p}}{\mathbf{l}}\otimes \mathbf{I_m}).\label{Pp}
\end{eqnarray}
Denote
$\mathbf{\hat{A}_1}={\rm Diag}(\hat{\mathcal{A}}_1),\mathbf{\hat{A}_2}={\rm Diag}(\hat{\mathcal{A}}_2), \mathbf{\hat{B}_1}={\rm Diag}(\hat{\mathcal{B}}_1),\;\mbox{\rm and}\;\mathbf{\hat{B}_2}={\rm Diag}(\hat{\mathcal{B}}_2).$
Then, from Proposition \ref{pro:associative} (i), (v), (vi) and \eqref{Pp}, it follows
\begin{equation}\label{part1}
\begin{array}{lcl}
&&[(-\mathbf{F_{p}}\mathbf{l}\otimes \mathbf{I_{m}})\;{\rm Diag}(\hat{\mathcal{A}}_1)\;(\mathbf{F_{p}}\mathbf{l}\otimes\mathbf{I_{n}})]\cdot
[(-\mathbf{F_{p}}\mathbf{l}\otimes \mathbf{I_{n}})\;{\rm Diag}(\hat{\mathcal{B}}_1)\;(\mathbf{F_{p}}\mathbf{l}\otimes\mathbf{I_{s}})]\\
&=&(-\mathbf{F_{p}}\mathbf{l}\otimes \mathbf{I_{m}})\;\mathbf{\hat{A}_1}\;(\mathbf{F_{p}}\mathbf{l}\otimes\mathbf{I_{n}})
(-\mathbf{F_{p}}\mathbf{l}\otimes \mathbf{I_{n}})\;\mathbf{\hat{B}_1}\;(\mathbf{F_{p}}\mathbf{l}\otimes\mathbf{I_{s}})\\
&=&(-\mathbf{F_{p}}\mathbf{l}\otimes \mathbf{I_{m}})\;\mathbf{\hat{A}_1}\;[(\mathbf{F_{p}}\mathbf{l}\otimes\mathbf{I_{n}})\;
(-\mathbf{F_{p}}\mathbf{l}\otimes \mathbf{I_{n}})]\;\mathbf{\hat{B}_1}\;(\mathbf{F_{p}}\mathbf{l}\otimes\mathbf{I_{s}})\\
&=&(-\mathbf{F_{p}}\mathbf{l}\otimes \mathbf{I_{m}})\;\mathbf{\hat{A}_1}\;\mathbf{\hat{B}_1}\;(\mathbf{F_{p}}\mathbf{l}\otimes\mathbf{I_{s}})\\
&=&[(-\mathbf{F^*_{p}}\mathbf{l}\otimes \mathbf{I_{m}})\;(\mathbf{{P}_{p}}\otimes \mathbf{I_{m}})]\cdot\mathbf{\hat{A}_1}\;\mathbf{\hat{B}_1}
\cdot[(\mathbf{{P}_{p}}\otimes\mathbf{I_{s}})\;(\mathbf{F^*_{p}}\mathbf{l}\otimes\mathbf{I_{s}})]\\
&=&(-\mathbf{F^*_{p}}\mathbf{l}\otimes \mathbf{I_{m}})\cdot[(\mathbf{{P}_{p}}\otimes \mathbf{I_{m}})\;\mathbf{\hat{A}_1}\;\mathbf{\hat{B}_1}\;(\mathbf{{P}_{p}}\otimes \mathbf{I_{s}})]\cdot(\mathbf{F^*_{p}}\mathbf{l}\otimes\mathbf{I_{s}})\\
&=&(-\mathbf{F^*_{p}}\mathbf{l}\otimes \mathbf{I_{m}})\cdot[(\mathbf{{P}_{p}}\otimes \mathbf{I_{m}})\;{\rm Diag}({\hat{\mathcal{A}}_1})\;{\rm Diag}({\hat{\mathcal{B}}_1})\;(\mathbf{{P}_{p}}\otimes\mathbf{I_{s}})]\cdot
(\mathbf{F^*_{p}}\mathbf{l}\otimes\mathbf{I_{s}}).
\end{array}
\end{equation}
From Proposition \ref{pro:associative} (ii), (v), (vi) and \eqref{Pp}, it follows
\begin{equation}\label{part2}
\begin{array}{ll}
&[(-\mathbf{F^*_{p}}\mathbf{l}\otimes \mathbf{I_{m}})\;({\rm Diag}(\hat{\mathcal{A}}_2)\mathbf{j})\;(\mathbf{F^*_{p}}\mathbf{l}\otimes\mathbf{I_{n}})]\cdot
[(-\mathbf{F^*_{p}}\mathbf{l}\otimes \mathbf{I_{n}})\;({\rm Diag}(\hat{\mathcal{B}}_2)\mathbf{j})\;(\mathbf{F^*_{p}}\mathbf{l}\otimes\mathbf{I_{s}})]\;\\
=&[(-\mathbf{F^*_{p}}\mathbf{l}\otimes \mathbf{I_{m}})\;(\mathbf{\mathbf{\hat{A}_2}}\mathbf{j})\;(\mathbf{F^*_{p}}\mathbf{l}\otimes\mathbf{I_{n}})]\cdot
\{(-\mathbf{F_{p}}\mathbf{l}\otimes \mathbf{I_{n}})\;[(\mathbf{{P}_{p}}\otimes \mathbf{I_{n}})\;(\mathbf{\hat{B}_2}\mathbf{j})\;(\mathbf{{P}_{p}}\otimes \mathbf{I_{s}})]\;(\mathbf{F_{p}}\mathbf{l}\otimes\mathbf{I_{s}})\}\;\\
=&[(-\mathbf{F^*_{p}}\mathbf{l}\otimes \mathbf{I_{m}})\;(\mathbf{\hat{A}_2}\mathbf{j})\;(\mathbf{F^*_{p}}\mathbf{l}\otimes\mathbf{I_{n}})]\cdot (-\mathbf{F^*_{p}}\mathbf{l}\otimes \mathbf{I_{n}})\cdot[(\mathbf{{P}_{p}}\otimes \mathbf{I_{n}})\;(\mathbf{\hat{B}_2}\mathbf{j})\;(\mathbf{{P}_{p}}\otimes \mathbf{I_{s}})]\cdot(\mathbf{F^*_{p}}\mathbf{l}\otimes\mathbf{I_{s}})\\
=&(-\mathbf{F^*_{p}}\mathbf{l}\otimes \mathbf{I_{m}})\;(\mathbf{\mathbf{\hat{A}_2}}\mathbf{j})\cdot[(\mathbf{F_{p}}\mathbf{l}\otimes\mathbf{I_{n}})
\;(-\mathbf{F_{p}}\mathbf{l}\otimes \mathbf{I_{n}})]\cdot[(\mathbf{{P}_{p}}\otimes \mathbf{I_{n}})\;(\mathbf{\hat{B}_2}\mathbf{j})\;(\mathbf{{P}_{p}}\otimes \mathbf{I_{s}})]\cdot(\mathbf{F^*_{p}}\mathbf{l}\otimes\mathbf{I_{s}})\\
=&(-\mathbf{F^*_{p}}\mathbf{l}\otimes \mathbf{I_{m}})\cdot[-\overline{\mathbf{\hat{A}_2}}\;(\mathbf{{P}_{p}}\otimes \mathbf{I_{n}})\;\hat{\mathbf{B}}_2\;(\mathbf{{P}_{p}}\otimes \mathbf{I_{s}})]\cdot(\mathbf{F^*_{p}}\mathbf{l}\otimes\mathbf{I_{s}})\\
=&(-\mathbf{F^*_{p}}\mathbf{l}\otimes \mathbf{I_{m}})\cdot[-{\rm Diag}(\overline{{\hat{\mathcal{A}}_2}})\;(\mathbf{{P}_{p}}\otimes \mathbf{I_{n}})\;{\rm Diag}(\hat{\mathcal{B}}_2)\;(\mathbf{{P}_{p}}\otimes\mathbf{I_{s}})]\cdot
(\mathbf{F^*_{p}}\mathbf{l}\otimes\mathbf{I_{s}}).
\end{array}
\end{equation}
From Proposition \ref{pro:associative} (i), (v), (vi) and \eqref{Pp}, it follows
\begin{equation}\label{part3}
\begin{array}{lcl}
&&[(-\mathbf{F_{p}}\mathbf{l}\otimes \mathbf{I_{m}})\;{\rm Diag}(\hat{\mathcal{A}}_1)\;(\mathbf{F_{p}}\mathbf{l}\otimes\mathbf{I_{n}})]\cdot
[(-\mathbf{F^*_{p}}\mathbf{l}\otimes \mathbf{I_{n}})\;({\rm Diag}(\hat{\mathcal{B}}_2)\mathbf{j})\;(\mathbf{F^*_{p}}\mathbf{l}\otimes\mathbf{I_{s}})]\\
&=&\{(-\mathbf{F^*_{p}}\mathbf{l}\otimes \mathbf{I_{m}})\cdot[(\mathbf{{P}_{p}}\otimes \mathbf{I_{m}})\;\mathbf{\hat{A}_1}\;(\mathbf{{P}_{p}}\otimes \mathbf{I_{n}})]\cdot(\mathbf{F^*_{p}}\mathbf{l}\otimes\mathbf{I_{n}})\}\cdot
[(-\mathbf{F^*_{p}}\mathbf{l}\otimes \mathbf{I_{n}})\;(\mathbf{\hat{B}_2}\mathbf{j})\;(\mathbf{F^*_{p}}\mathbf{l}\otimes\mathbf{I_{s}})]\;\\
&=&(-\mathbf{F^*_{p}}\mathbf{l}\otimes \mathbf{I_{m}})\cdot[(\mathbf{{P}_{p}}\otimes \mathbf{I_{m}})\;\mathbf{\hat{A}_1}\;(\mathbf{{P}_{p}}\otimes \mathbf{I_{n}})]\cdot[(\mathbf{F^*_{p}}\mathbf{l}\otimes\mathbf{I_{n}})\;
(-\mathbf{F^*_{p}}\mathbf{l}\otimes \mathbf{I_{n}})]\cdot(\mathbf{\hat{B}_2}\mathbf{j})\cdot(\mathbf{F^*_{p}}\mathbf{l}\otimes\mathbf{I_{s}})\\
&=&(-\mathbf{F^*_{p}}\mathbf{l}\otimes \mathbf{I_{m}})\cdot[(\mathbf{{P}_{p}}\otimes \mathbf{I_{m}})\;\overline{\mathbf{\hat{A}_1}}\;(\mathbf{{P}_{p}}\otimes \mathbf{I_{n}})\;(\mathbf{\hat{B}_2}\mathbf{j})]\cdot(\mathbf{F^*_{p}}\mathbf{l}\otimes\mathbf{I_{s}})\\
&=&(-\mathbf{F^*_{p}}\mathbf{l}\otimes \mathbf{I_{m}})\cdot[(\mathbf{{P}_{p}}\otimes \mathbf{I_{m}})\;{\rm Diag}(\overline{\hat{\mathcal{A}}_1})\;(\mathbf{{P}_{p}}\otimes \mathbf{I_{n}})\;({\rm Diag}(\hat{\mathcal{B}}_2)\mathbf{j})]\cdot(\mathbf{F^*_{p}}\mathbf{l}\otimes\mathbf{I_{s}}).
\end{array}
\end{equation}
From Proposition \ref{pro:associative} (ii), (v) and (vi), it follows
\begin{equation}\label{part4}
\begin{array}{lcl}
&&[(-\mathbf{F^*_{p}}\mathbf{l}\otimes \mathbf{I_{m}})\;({\rm Diag}(\hat{\mathcal{A}}_2)\mathbf{j})\;(\mathbf{F^*_{p}}\mathbf{l}\otimes\mathbf{I_{n}})]\cdot
[(-\mathbf{F_{p}}\mathbf{l}\otimes \mathbf{I_{n}})\;{\rm Diag}(\hat{\mathcal{B}}_1)\;(\mathbf{F_{p}}\mathbf{l}\otimes\mathbf{I_{s}})]\\
&=&[(-\mathbf{F^*_{p}}\mathbf{l}\otimes \mathbf{I_{m}})\;({\rm Diag}(\hat{\mathcal{A}}_2)\mathbf{j})\;(\mathbf{F^*_{p}}\mathbf{l}\otimes\mathbf{I_{n}})]\cdot
(-\mathbf{F^*_{p}}\mathbf{l}\otimes \mathbf{I_{n}})\cdot{\rm Diag}(\hat{\mathcal{B}}_1)\cdot(\mathbf{F^*_{p}}\mathbf{l}\otimes\mathbf{I_{s}})\\
&=&(-\mathbf{F^*_{p}}\mathbf{l}\otimes \mathbf{I_{m}})\cdot({\rm Diag}(\hat{\mathcal{A}}_2)\mathbf{j})\cdot[(\mathbf{F_{p}}\mathbf{l}\otimes\mathbf{I_{n}})\;
(-\mathbf{F_{p}}\mathbf{l}\otimes \mathbf{I_{n}})]\cdot{\rm Diag}(\hat{\mathcal{B}}_1)\cdot(\mathbf{F^*_{p}}\mathbf{l}\otimes\mathbf{I_{s}})\\
&=&(-\mathbf{F^*_{p}}\mathbf{l}\otimes \mathbf{I_{m}})\cdot[{\rm Diag}(\hat{\mathcal{A}}_2){\rm Diag}(\hat{\mathcal{B}}_1)\mathbf{j}]\cdot(\mathbf{F^*_{p}}\mathbf{l}\otimes\mathbf{I_{s}}).
\end{array}
\end{equation}
Thus, by \eqref{bcirc(C)}, \eqref{part1}, \eqref{part2}, \eqref{part3}, \eqref{part4}, \eqref{eq:A} and \eqref{eq:B}, we obtain that
$$\begin{array}{c}
bcirc(\mathcal{C})=(\mathbf{F^*_{p}}\mathbf{l}\otimes \mathbf{I_{m}})\;[{\rm Diag}({\hat{\mathcal{C}}_1})+{\rm Diag}({\hat{\mathcal{C}}_2})\mathbf{j}]\;(-\mathbf{F^*_{p}}\mathbf{l}\otimes\mathbf{I_{s}}),
\end{array}$$
i.e., \eqref{eq:T-product} holds.\ep

\begin{remark}
It is easy to see that \eqref{C} and \eqref{C2} are equivalent to
\begin{eqnarray}
{\hat{\mathcal{C}}}_1(:,:,i)&=& {\hat{\mathcal{A}}_1}(:,:,p+2-i)\;\hat{\mathcal{B}}_1(:,:,p+2-i)
-\overline{{\hat{\mathcal{A}}_2}}(:,:,i)\;\hat{\mathcal{B}}_2(:,:,p+2-i),\label{TC} \\
{\hat{\mathcal{C}}_2}(:,:,i)&=&\overline{\hat{\mathcal{A}}_1}(:,:,p+2-i)\;
\hat{\mathcal{B}}_2(:,:,i)+\hat{\mathcal{A}}_2(:,:,i)\;\hat{\mathcal{B}}_1(:,:,i),\label{TC2}
\end{eqnarray}
where $\hat{\mathcal{A}}(:,:,p+1):=\hat{\mathcal{A}}(:,:,1)$ and  $\hat{\mathcal{B}}(:,:,p+1):=\hat{\mathcal{B}}(:,:,1)$.

Let $\hat{\mathcal{C}}=\hat{\mathcal{C}}_1+\hat{\mathcal{C}}_2\mathbf{j}.$
Then from \eqref{eq:T-product} and Remark
\ref{rem-2}, we can obtained that $\mathcal{C}$ can be computed by applying the inverse fast Fourier transform to $\overline{\hat{\mathcal{C}}}=\overline{\hat{\mathcal{C}}_1}-\hat{\mathcal{C}}_2\mathbf{j}$ along the third dimension, i.e., implementing \begin{equation}\label{eq:ifft-C}
\mathcal{C}=\operatorname{ifft}(\overline{{\hat{\mathcal{C}}_1}},[\;]\;,3)-
\operatorname{ifft}(\hat{\mathcal{C}}_2,[\;]\;,3)\mathbf{j}
:=\operatorname{ifft}(\overline{{\hat{\mathcal{C}}}},[\;]\;,3).
\end{equation}
\end{remark}

With Theorem \ref{theo:T-product}, we give a way to calculate the T-product between two quaternion tensors via FFTs, which is shown in Algorithm 5.1.

\begin{center}
{\begin{tabular}{rl}
  \hline\noalign{\smallskip}
{\bf Algorithm}&{\bf 5.1: Computing the T-product of quaternion tensors via FFTs} \\
\noalign{\smallskip}\hline\noalign{\smallskip}
{\bf \;\;Input:}& $\mathcal{A}\in \mathbb{Q}^{m\times n\times p}$ and $\mathcal{B}\in \mathbb{Q}^{n\times s \times p}$.\\
{\;\;\bf Output:}& $\mathcal{C}\in \mathbb{Q}^{m\times s\times p}$: $\mathcal{C}=\mathcal{A}\ast\mathcal{B}$.\\
{\;\;\bf Step 1:}&Compute $\hat{\mathcal{A}}=\operatorname{fft}(\overline{\mathcal{A}},[\;]\;,3)$ and $\hat{\mathcal{B}}=\operatorname{fft}(\overline{\mathcal{B}},[\;]\;,3)$ as in {\eqref{eq:fftA-B}}.\\
{\;\;\bf Step 2:}& Let $\hat{\mathcal{A}}(:,:,p+1):=\hat{\mathcal{A}}(:,:,1)$ and  $\hat{\mathcal{B}}(:,:,p+1):=\hat{\mathcal{B}}(:,:,1)$.\\
& For $i=1:p,$ compute \\
&\quad\quad$\hat{\mathcal{C}}_1(:,:,i)$ and $\hat{\mathcal{C}}_2(:,:,i)$ as in \eqref{TC} and \eqref{TC2}, respectively,\\
& end for.\\
{\;\;\bf Step 3:}& Compute $\mathcal{C}=\operatorname{ifft}(\overline{\hat{\mathcal{C}}},[\;]\;,3)$ as in \eqref{eq:ifft-C}.\\
\noalign{\smallskip}  \hline
\end{tabular}}
\end{center}

%
In the following, we show that the computational magnitude of calculating $\mathcal{A}\ast\mathcal{B}$ by Algorithm 5.1 is $\frac{1}{p}$ of the one of direct calculating $\mathcal{A}\ast\mathcal{B}$ by Definition \ref{T-product}, which can be seen from the proof of Theorem \ref{theo:computational-magnitude} and  calculations of  Example \ref{exa:2}.

\begin{theorem}\label{theo:computational-magnitude}
Let $\mathcal{A}\in \mathbb{Q}^{m\times n\times p}$ and $\mathcal{B}\in \mathbb{Q}^{n\times s\times p}$. Then the computational magnitude of calculating $\mathcal{A}\ast\mathcal{B}$ by Algorithm 5.1 is $\mathbf{O}(mnsp)$, and the computational magnitude of direct calculating $\mathcal{A}\ast\mathcal{B}$ by Definition \ref{T-product} is $\mathbf{O}(mnsp^2)$.
\end{theorem}

{\bf Proof.} First, we analyze the computational magnitude of calculating $\mathcal{A}\ast\mathcal{B}$ by Algorithm 5.1.

For the computation of the fft and ifft in {\bf Step 1} and {\bf Step 3},
as shown in \cite{GV-2013,KM-2011}, the computational magnitudes of fft of an $m\times n\times p$ complex tensor is $\mathbf{O}(mnp\log p)$. Noting that the fft of a quaternion tensor requires two times ffts of the corresponding complex tensor, thus the computational magnitude of fft of $\overline{\mathcal{A}}$ and fft of $\overline{\mathcal{B}}$ in {\bf Step 1} of Algorithm 5.1 are $\mathbf{O}(2mnp\log p)$ and $\mathbf{O}(2nsp\log p)$, respectively. Similarly, the computational magnitude of ifft of $\overline{\hat{\mathcal{C}}}$ in {\bf Step 3} of Algorithm 5.1 is $\mathbf{O}(2msp\log p)$.

For the computation of the multiplication of complex matrices in {\bf Step 2},
the main cost of {\bf Step 2} in Algorithm 5.1 is to compute multiplications of complex matrices ${\rm Diag}(\hat{\mathcal{A}}_1)$ and ${\rm Diag}(\hat{\mathcal{B}}_1)$, ${\rm Diag}(\overline{{\hat{\mathcal{A}}_2}})$ and $[(\mathbf{{P}_{p}}\otimes \mathbf{I_{n}}){\rm Diag}(\hat{\mathcal{B}}_2)(\mathbf{{P}_{p}}\otimes\mathbf{I_{s}})]\;$, ${\rm Diag}(\hat{\mathcal{B}}_2)$ and $[(\mathbf{{P}_{p}}\otimes \mathbf{I_{m}}){\rm Diag}(\overline{\hat{\mathcal{A}}_1})(\mathbf{{P}_{p}}\otimes\mathbf{I_{n}})]\;$,
${\rm Diag}(\hat{\mathcal{A}}_2)$ and ${\rm Diag}(\hat{\mathcal{B}}_1)$. Note that ${\rm Diag}(\hat{\mathcal{A}}_1)$, ${\rm Diag}(\hat{\mathcal{A}}_2)$, ${\rm Diag}(\overline{\hat{\mathcal{A}}_1})$, ${\rm Diag}(\overline{\hat{\mathcal{A}}_2})$ are block diagonal matrices with $p$ blocks of $m\times n$, and ${\rm Diag}(\hat{\mathcal{B}}_1)$, ${\rm Diag}(\hat{\mathcal{B}}_2)$, ${\rm Diag}(\overline{\hat{\mathcal{B}}_1})$, ${\rm Diag}(\overline{\hat{\mathcal{B}}_2})$ are block diagonal matrices with $p$ blocks of $n\times s$, and $[(\mathbf{{P}_{p}}\otimes \mathbf{I_{n}}){\rm Diag}(\hat{\mathcal{B}}_2)(\mathbf{{P}_{p}}\otimes\mathbf{I_{s}})]\;$ ($[(\mathbf{{P}_{p}}\otimes \mathbf{I_{m}}){\rm Diag}(\overline{\hat{\mathcal{A}}_1})(\mathbf{{P}_{p}}\otimes\mathbf{I_{n}})]\;$) is achieved by permutating the $r$-th block and the $(p+2-r)$-th block for all $r=2,3,\ldots,p$. Thus, the calculations in {\bf Step 2} of Algorithm 5.1 requires $4p$ multiplications between an $m\times n$ complex matrix and an $n\times s$ complex matrix. Noting that one multiplication between two complex numbers requires $4$ real multiplications and $2$ real additions, and one addition between two complex numbers requires $2$ real additions, thus the computation of {\bf Step 2} of Algorithm 5.1 requires $16mnsp$ real multiplications and $8m(n-1)sp+8mnsp$ real additions in total. Hence, the computational magnitude of {\bf Step 2} of Algorithm 5.1 is $\mathbf{O}(16mnsp)$.

So, the computational magnitude of calculating $\mathcal{A}\ast\mathcal{B}$ by Algorithm 5.1 is $\mathbf{O}(mnsp)$.

Second, we show the computational magnitude of direct calculating $\mathcal{A}\ast\mathcal{B}$ by its definition.
The computational magnitude of direct calculating $\mathcal{A}\ast\mathcal{B}$ is obtained by multiplying the $mp\times np$ circulant matrix $bcirc(\mathcal{A})$ and the $np\times s$ block column vector $unfold(\mathcal{B})$. Thus, direct calculating $\mathcal{A}\ast\mathcal{B}$ requires $mnsp^2$ multiplications and $mps(np-1)$ additions between quaternions. Noting that one multiplication between two quaternions requires $16$ real multiplications and $12$ real additions, and one addition between two quaternions requires $4$ real additions. Hence, the T-product between an $m\times n\times p$ quaternion tensor and an $n\times s\times p$ quaternion tensor requires $16mnsp^2$ real multiplications and $4mps(np-1)+12mnsp^2$ real additions in total.

So, the computational magnitude of direct calculating $\mathcal{A}\ast\mathcal{B}$ by Definition \ref{T-product} is $\mathbf{O}(mnsp^2)$. \ep

From the complexity analysis in Theorem \ref{theo:computational-magnitude}, it can be seen that Algorithm 5.1 can speed up almost $p$ times of the calculation of $\mathcal{A}\ast\mathcal{B}$ with $p$ frontal slices, which is verified by Example \ref{exa:2}.

\begin{example}\label{exa:2}
Let $\mathcal{A}\in \mathbb{Q}^{m\times n\times p}$ and
$\mathcal{B}\in \mathbb{Q}^{n\times s\times p}$ be random generated third-order quaternion tensors.
We test the effectiveness of Algorithm 5.1 to calculate $\mathcal{A}\ast\mathcal{B}$. The numerical results are shown in Table \ref{tab:T-product}, where in each pair ``$(m,n,s)$'', $m\times n$ is the size of each frontal slice of $\mathcal{A}$ and $n\times s$ is the size of each frontal slice of $\mathcal{B}$; ``$p$'' is the number of frontal slices of $\mathcal{A}$ ($\mathcal{B}$); ``Time\_Algo\_5.1'' means the time in seconds spent for calculating $\mathcal{A}\ast\mathcal{B}$ by applying Algorithm 5.1; ``Time\_Direct'' means the time in seconds spent for calculating $\mathcal{A}\ast\mathcal{B}$ directly according to its definition; ``error'' means the relative error between calculation results of $\mathcal{A}\ast\mathcal{B}$ by applying Algorithm 5.1 and directly calculating by Definition \ref{T-product}. The computation is performed on a HP Laptop with CPU of 2.6 GHz and RAM of 8.0 GB.

\begin{table}[!htbp]\;
  \centering{
     \begin{tabular}{ccccccccc}
     \hline\noalign{\smallskip}
     ($m,n,s$)    & $p$      & Time\_Algo\_5.1(s)  &  Time\_Direct(s)   & error   \\
       \noalign{\smallskip}   \hline\noalign{\smallskip}
     (50,50,50)  & 10     &  0.0197                 & 0.0623             &$5.7413e^{-17}$        \\
     (50,50,50)  & 50     &  0.0531                 & 1.1554             &$1.0763e^{-16}$        \\
     (50,50,50)  & 100    &  0.1064                 & 4.7656             &$1.4844e^{-16}$       \\
     {\bf(50,50,50)}  & {\bf 200}    &  {\bf 0.3801}  & {\bf71.4014}     &${\bf 2.0711e^{-16}}$     \\ (100,100,100)  & 10  &  0.0609                 & 0.3176             &$4.8256e^{-17}$        \\
     (100,100,100)  & 50  &  0.2928                 & 8.6356             &$4.2408e^{-17}$       \\
     {\bf(100,100,100)}  & {\bf 100} &  {\bf 0.8533}  & {\bf78.9583}     &${\bf 4.3087e^{-17}}$      \\
     \hline
\end{tabular}}
\caption{Comparisons of computing $\mathcal{A}\ast\mathcal{B}$ by Algorithm 5.1 and by Definition \ref{T-product}}
\label{tab:T-product}
\end{table}

It can be seen from Table \ref{tab:T-product}:
\begin{itemize}
\item for the calculation of $\mathcal{A}\ast\mathcal{B}$ with $\mathcal{A},\mathcal{B}\in \mathbb{Q}^{50\times 50\times 200}$, applying Algorithm 5.1 only takes {\bf 0.3801} seconds, almost ${\bf\frac{1}{200}}$ of the time taken by directly calculating by Definition \ref{T-product}, which takes {\bf 71.4014} seconds, and the relative error between two results reaches ${\bf 10^{-16}}$;
\item for the calculation of $\mathcal{A}\ast\mathcal{B}$ with $\mathcal{A},\mathcal{B}\in \mathbb{Q}^{100\times 100\times 100}$, applying Algorithm 5.1 only takes {\bf 0.8533} seconds, almost ${\bf \frac{1}{100}}$ of the time taken by directly calculating by Definition \ref{T-product}, which takes {\bf 78.9583} seconds, and the relative error between two results reaches ${\bf 10^{-17}}$.
\end{itemize}
\end{example}

%

\section{Concluding remarks}
Since T-product based third-order real tensors methods based on the block diagonalization of block circulant real matrices via FFTs have achieved great success in the color image and gray video processing, our goal was to achieve the block diagonalization of block circulant quaternion matrices with the calculation convenience of FFTs in this paper. We first studied the block diagonalization of block circulant quaternion matrices in the quaternion domain. We showed that a circulant quaternion matrix of size $p\times p$ cannot be diagonalized by the unitary DFT matrix $\mathbf{F_p}$, nor by the unitary quaternion matrices $\mathbf{F_p}\mathbf{j}$ and $(\mathbf{F_p}+\mathbf{F_p}\mathbf{j})/\sqrt{2}$. After that, we proved sufficient and necessary conditions for a unitary quaternion matrix being a diagonalization matrix of circulant quaternion matrices. These studies in the quaternion domain showed that it was a hard work to find unitary quaternion matrices to achieve the block diagonalization of block circulant quaternion matrices in the quaternion domain, and FFTs cannot be applied to obtain result block diagonal matrices of block circulant quaternion matrices in the quaternion domain even if such unitary quaternion matrices were found. Owing to that, we then studied the diagonalization of circulant quaternion matrices in the octonion domain. We proved that a circulant (respectively, block circulant) quaternion matrix of size $p\times p$ (respectively, of size $mp\times mp$) can be transformed into a diagonal (respectively, block diagonal) quaternion matrix by unitary octonion matrices, at a cost of $\mathbf{O}(p\log p)$ (respectively, $\mathbf{O}(mnp\log p)$) via FFTs. Further, an algorithm to fast calculate the T-product of third-order quaternion tensors via FFTs was proposed by establishing the relationship between the T-product calculation and the block diagonalization of block circulant matrices. Numerical calculations demonstrated that the proposed method speeds up almost $p$ times for computing the T-product between third-order quaternion tensors with $p$ frontal slices.

With the block diagonalization of block circulant quaternion matrices achieved in the octonion domain and the fast calculation of T-product between third-order quaternion tensors, T-product based third-order quaternion tensors methods are possible now. Then how to extend these efficent T-product based third-order real tensors methods \cite{KMP-2008,B-2010,HK-2010,KM-2011,KBHH-2013,ZEAHK-2014,ZA-2017,SNZ-2020,ZBN-2020,
Lu-2019,Lu-2016} to the quaternion case to apply for various problems, such as the color videos processing? It is worthwhile to further study.

Besides, it is known that a lot of matrix theory that holds for complex matrices cannot be extended to the quaternion case, due to the non-commutativity of quaternion algebra. But analogue to the discussion in Section \ref{Sect.DQCM}, with the help of the special rule of the multiplication between special octonions and quaternions as shown in \eqref{eq:relation-remark}, it is not difficult to deduce that for any $\mathbf{P}$ and $\mathbf{Q}$ being complex matrices and $\mathbf{A}$ being a quaternion matrix:
$$(\mathbf{P}\mathbf{l})\;\mathbf{A}\;(\mathbf{Q}\mathbf{l})=
(\mathbf{P}\mathbf{jl})\;\mathbf{A}\;(\mathbf{Q}\mathbf{jl})
=[\mathbf{P}(\mathbf{l}+\mathbf{jl})/\sqrt{2}]\;\mathbf{A}\;[\mathbf{Q}(\mathbf{l}+\mathbf{jl})/\sqrt{2}]
=-\mathbf{P}\;\overline{\mathbf{A_1}}\;\overline{\mathbf{Q}}+(\mathbf{P}\;\mathbf{A_2}\;
\overline{\mathbf{Q}})\;\mathbf{j},$$
which means that if complex matrices $\mathbf{P}$ and $\overline{\mathbf{Q}}$ can transform special structured complex matrices  into diagonal complex matrices, then octonion matrices $\mathbf{P}\mathbf{l}$ (or $\mathbf{P}\mathbf{jl}$, $\mathbf{P}(\mathbf{l}+\mathbf{jl})/\sqrt{2}$, respectively) and $\mathbf{Q}\mathbf{l}$ (or $\mathbf{Q}\mathbf{jl}$, $\mathbf{Q}(\mathbf{l}+\mathbf{jl})/\sqrt{2}$, respectively) can transform quaternion matrices with same special structures into diagonal quaternion matrices. This encourages us to raise an issue: whether can more matrix theory be extended from the complex domain to the quaternion domain through the participation of octonion matrices? It is also deserve to further study.



\end{document}